\documentclass[a4paper,10pt]{article}

\usepackage[caption=false]{subfig}
\usepackage[utf8]{inputenc}
\usepackage[T1]{fontenc}
\usepackage[english]{babel}
\usepackage{lmodern}
\usepackage{amsmath,amsfonts,dsfont,amssymb}
\usepackage{array}
\usepackage{pifont}
\usepackage{csquotes}
\usepackage{ragged2e}
\usepackage{authblk}
\usepackage{graphicx}
\usepackage{algorithm}
\usepackage[noend]{algorithmic}
\usepackage[numbers]{natbib}
\usepackage{fullpage}
\usepackage{xcolor}

\DeclareMathOperator{\Var}{Var}
\DeclareMathOperator{\cov}{Cov}

\DeclareMathOperator{\HSIC}{HSIC}

\DeclareMathOperator{\R2}{R}
\DeclareMathOperator{\HS}{HS}

\DeclareMathOperator{\Pval}{P}
\DeclareMathOperator{\P1}{P_1}
\DeclareMathOperator{\Pk}{P_k}
\DeclareMathOperator{\Pd}{P_d}
\DeclareMathOperator{\pval}{p}
\DeclareMathOperator{\p1}{p_1}
\DeclareMathOperator{\pk}{p_k}
\DeclareMathOperator{\pd}{p_d}
\DeclareMathOperator{\F}{F}
\DeclareMathOperator{\Gk}{G_k}

\DeclareMathOperator{\Res}{\mathcal{R}}
\DeclareMathOperator{\D1}{\mathcal{D}_1}
\DeclareMathOperator{\Dk}{\mathcal{D}_k}
\DeclareMathOperator{\Dd}{\mathcal{D}_d}
\DeclareMathOperator{\MMD}{MMD}
\DeclareMathOperator{\SuppZ}{Supp(Z)}
\DeclareMathOperator{\T}{T}
\DeclareMathOperator{\Tr}{Tr}
\DeclareMathOperator{\SKL}{SKL}
\DeclareMathOperator{\KL}{\mathbf{KL}}

\newcommand{\W}{W}

\newcommand{\tL}{\widetilde{L}}

\newcommand{\K}{K}
\newcommand{\M}{M}
\newcommand{\WW}{W}

\newtheorem{prop}{Proposition}

\newtheorem{remark}{Remark}

\usepackage[hyperfootnotes=false]{hyperref}
\hypersetup{
  colorlinks,
  citecolor=blue,
  linkcolor=red,
  urlcolor=gray}

\begin{document}
\title{Second-level global sensitivity analysis of numerical simulators with application to an accident scenario in a sodium-cooled fast reactor}

\author[1,2,3]{Anouar Meynaoui \footnote{anouar.meynaoui@gmail.com}}
\author[1,2]{Amandine Marrel}
\author[2,3]{Béatrice Laurent}

\affil[1]{CEA, DES, IRESNE, DER, Cadarache F-13108 Saint-Paul-Lez-Durance, France.}
\affil[2]{Institut de Mathématiques de Toulouse, UMR5219, F-31062 Toulouse, France.}
\affil[3]{Université de Toulouse; CNRS, INSA, F-31077 Toulouse, France.}

\maketitle

\abstract{Numerical simulators are widely used to model physical phenomena and global sensitivity analysis (GSA) aims at studying the global impact of the input uncertainties on the simulator output. To perform GSA, statistical tools based on inputs/output dependence measures are commonly used. We focus here on the Hilbert-Schmidt independence criterion (HSIC). Sometimes, the probability distributions modeling the uncertainty of inputs may be themselves uncertain and it is important to quantify their impact on GSA results. We call it here the second-level global sensitivity analysis (GSA2). However, GSA2, when performed with a Monte Carlo double-loop, requires a large number of model evaluations, which is intractable with CPU time expensive simulators. To cope with this limitation, we propose a new statistical methodology based on a Monte Carlo single-loop with a limited calculation budget. First, we build a unique sample of inputs and simulator outputs, from a well-chosen probability distribution of inputs. From this sample, we perform GSA for various assumed probability distributions of inputs by using weighted HSIC measures estimators. Statistical properties of these weighted estimators are demonstrated. Subsequently, we define 2$^{\text{nd}}$-level HSIC-based measures between the distributions of inputs and GSA results, which constitute GSA2 indices. The efficiency of our GSA2 methodology is illustrated on an analytical example, thereby comparing several technical options. Finally, an application to a test case simulating a severe accidental scenario on nuclear reactor is provided.}


\maketitle

\section{Introduction}
\label{Intro}

Numerical simulators are fundamental tools for understanding, modeling and predicting phenomena. They are widely used nowadays in several fields such as physics, chemistry and biology, but also in economics and social science. These numerical simulators take a large number of input parameters more or less uncertain, characterizing the studied phenomenon. Consequently, the output which is provided by the numerical simulator is also uncertain. It is therefore important to consider not only the nominal values of inputs, but also the set of all possible values in the range of variation of each uncertain input \cite{Derocquigny2008, Helton2006}. In the framework of a probabilistic approach, the inputs and the output are considered as random variables and their uncertainties are modeled by probability distributions. The objective is then to evaluate the impact of the input uncertainties on the variability of the output. For this, sensitivity analysis studies can be performed, using statistical methods based on a sample of realizations from the simulator. To choose these numerical simulations, experimental design techniques can be used~\cite{damblin2013numerical}.

\smallskip
\textbf{Generalities on sensitivity analysis.} Sensitivity analysis \cite{saltelli2008global} aims at determining how the variability of inputs contributes, qualitatively or quantitatively, to the output variability. Sensitivity analysis can yield a screening of the inputs, which consists in separating them into two subgroups: those that mainly influence the output (most influential inputs) and those whose influence on the output can be neglected. More generally, sensitivity analysis can be divided into two main areas:

\begin{itemize}
\item local sensitivity analysis (LSA) which studies the output variability for a small input variation around nominal values (reference values);
\item global sensitivity analysis (GSA) which studies the impact of the input uncertainties on the output, considering the whole range of input variation.
\end{itemize}
 
We focus here on GSA and we call it in the following, first-level GSA, denoted GSA1. 

\smallskip
\textbf{Use of dependence measures for GSA1.} Among GSA1 tools \cite{iooss2015review}, one of the most popular methods used in industrial applications is based on a variance decomposition of the output \cite{sobol1993sensitivity}. The sensitivity indices thus obtained by this decomposition are called Sobol' indices. These indices have the advantage of being easily interpretable but are in practice very expensive in computing time (several tens of thousands of simulations required). More recently, tools based on dependence measures have been proposed for GSA1 purpose \cite{da2015global}. These measures aim at quantifying, from a probabilistic point of view, the dependence between the output random variable and the input random variables. Among these measures, we can mention the $f$-divergence of Csiszár which, for a given input, compares the distribution of the output and its distribution when this input is fixed, thanks to a function with specific properties \cite{csiszar1972class}. Always on the same principle, the distance correlation is an other dependence measure which compares the characteristic function of a couple of random input/output variables, with the product of the joint characteristic functions of the two variables \cite{szekely2007measuring}. Last but not least, the Hilbert-Schmidt independence criterion denoted HSIC \cite{gretton2005measuring}, generalizes the notion of covariance between two random variables and takes into account a very large spectrum of forms of dependence between variables. Initially developed by statisticians \cite{gretton2005measuring} to perform independence tests, these dependence measures offer the advantage of having a low cost of estimation (in practice a few hundred simulations against several tens of thousands for Sobol' indices) and their estimation for all inputs does not depend on the number of inputs. In addition, recent work proposed by \cite{de2016new} showed the efficiency of these measures to perform a screening of the input variables, from various HSIC-based statistical tests of significance. Finally, HSIC measures can easily be extended to non-vector inputs (functional, categorical, etc.). For all these reasons, we will focus here on HSIC measures for GSA1 of numerical simulators. 

\smallskip
\textbf{Second-level input uncertainties and GSA2.} In some cases, the probability distributions characterizing the uncertain inputs may themselves be uncertain. This uncertainty may be related to a divergence of expert opinion on the probability distribution assigned to each input or a lack of information to characterize this distribution. The modeling of this lack of knowledge on input laws can take many forms:

\begin{itemize}
\item the type of the input distribution is uncertain (uniform, triangular, normal, ...); 
\item the distribution is known but its parameters are uncertain (\textit{e.g.}, known normal distribution with unknown mean and variance, eventually estimated on data).
\end{itemize}

In both cases, the resulting uncertainties on the input laws are referred to here as \textit{second-level uncertainties}. As part of a probabilistic approach, these uncertainties can be modeled by a probability law on a set of possible probability laws of inputs or by a probability law on the parameters of a given input law (\textit{e.g.} Gaussian distribution with probability law on mean and/or variance). In any case, these 2$^{\text{nd}}$-level uncertainties can significantly change the GSA1 results performed by HSIC or any other dependence measure. In this framework, the main purpose of \textit{second-level GSA denoted GSA2} is to answer the following questions: \enquote{What impact do 2$^{\text{nd}}$-level uncertainties have on the GSA1 results?} and \enquote{What are the most influential ones and those whose influence is negligible?}. The GSA2 results and conclusion can then be used to prioritize the characterization efforts on the inputs whose uncertainties on probability laws have the greatest impact on GSA1 results. Note that, we assume here that the inputs are independent and continuous random variables with a probability density function, denoted here pdf. 

\smallskip
\textbf{Practical problems raised by GSA2.} In practice, the realization of GSA2 raises several issues and technical obstacles. First, it is necessary to characterize GSA1 results, i.e. to define a representative quantity of interest in order to compare the results obtained for different uncertain input pdf. Then, the impact of each uncertain input pdf on this quantity of interest has to be evaluated. For this, sensitivity indices measuring the dependence between GSA1 results and each input pdf have to be defined. We propose to call them \emph{2$^{\text{nd}}$-level} GSA indices. In order to estimate these measures, an approach based on a \enquote{Monte Carlo double-loop} could be considered. In the outer loop, a Monte Carlo sample of input pdfs is generated, while the inner loop aims at evaluating the GSA1 results associated to each pdf. For each pdf selected in the outer loop, the inner loop consists in generating a Monte Carlo sample of simulations (set of inputs/output) and to compute GSA1 results. The process is repeated for each input pdf. At the end of the outer loop, the impact of input pdf on the GSA1 results can be observed and quantify by computing 2$^{\text{nd}}$-level GSA. Unfortunately, this type of double-loop approach requires in practice a very large number of simulations which is intractable for time expensive computer simulators. Therefore, other less expensive approaches must be developed. 

\smallskip
To answer these different issues (choice of the quantity of interest, definition of 2$^{\text{nd}}$-level sensitivity indices and reduction of the budget of simulations), we propose in this paper a \enquote{single-loop} Monte Carlo methodology for GSA2 based on both 1$^{\text{st}}$-level and 2$^{\text{nd}}$-level HSIC dependence measures. Note that this work was initiated in the framework of Meynaoui's PhD \cite{meynaoui2019PhD}, the interested reader could find more technical elements and detailed demonstrations in this document.

\smallskip
The paper is organized as follows. In Section \ref{Foreword}, we introduce HSIC measures, before presenting the statistical estimators of these measures, as well as the associated characteristics (bias, variance and asymptotic law). Then, we show that these measures can be formulated and estimated with a sample generated from a different distribution than the \textit{prior} distribution of the inputs. For this, new estimators are proposed and their characteristics are detailed, these new estimators being a key point for the proposed GSA2 methodology. In Section \ref{NEW-GSA2}, the full methodology for GSA2 is presented: a single inputs/output sample is used, taking advantage of the new HSIC estimators. The GSA2 principle and the related practical issues are first introduced. The general algorithm is then detailed, followed by dedicated sections focusing on major technical elements. In Section \ref{Application}, the methodology is illustrated on an analytical example, thereby comparing different options and technical choices of the methodology. Finally, an application on a test case simulating a severe accidental scenario on a nuclear reactor is~proposed.

\section{Statistical inference around Hilbert-Schmidt dependence measures (HSIC)}
\label{Foreword}

Throughout the rest of this document, the numerical model is represented by the relation: 
\begin{equation*}
Y = \mathcal{M} \left( X_1 , \ldots , X_d \right),
\end{equation*}
where $X_1, \ldots ,X_d$ and $Y$ are respectively the $d$ uncertain inputs and the uncertain output, evolving in one-dimensional real sets respectively denoted $\mathcal {X}_1, \ldots, \mathcal{X}_d$ and $ \mathcal{Y}$. $\mathcal{M}$ denotes the numerical simulator. We note $\mathbf{X} = \left(X_1, \ldots, X_d \right)$ the vector of inputs. As part of the probabilistic approach, the $d$ inputs are considered as continuous and independent random variables with known densities. These densities are respectively denoted $f_1, \ldots , f_d$. Finally, $f : (x_1 , \ldots , x_d) \mapsto f_1 (x_1) \times \ldots \times f_d (x_d)$ denotes the density of the random vector $\mathbf{X}$. As the model $\mathcal{M}$ is not known analytically, a direct computation of the output probability density as well as dependence measures between $\mathbf{X}$ and $Y$ is impossible. Only observations (or realizations) of $\mathcal{M}$ are available. It is therefore assumed in the following that we have a $n$-sample of inputs and associated  outputs $\left(\mathbf{X}^{(i)} , Y^{(i)} \right)_{1 \leq i \leq n}$, where $Y^{(i)} = \mathcal{M} (\mathbf{X}^{(i)})$.

\subsection{Review on HSIC measures}
\label{HSIC}
After introducing their theoretical definition, the estimation of HSIC dependence measures and their use for GSA1 are detailed.   
 
\subsubsection{Definition and description}
\label{Principle-definition}
To define the HSIC measure between $X_k$ and $Y$, where $k \in \lbrace 1, \ldots , d \rbrace$, \cite{gretton2005measuring} associate to $X_k$ a \emph{reproducing kernel Hilbert space} (denoted RKHS, see \cite{aronszajn1950theory} for more details) $\mathcal{H}_k$ composed of functions mapping from $\mathcal {X}_k$ to $\mathbb{R}$ and characterized by a kernel $l_k$. The same transformation is carried out for $Y,$ considering a RKHS denoted $\mathcal{G}$ and a kernel $l$. The scalar products on $\mathcal{H}_k$ and $\mathcal{G}$ are respectively denoted $\langle . , . \rangle_{\mathcal{H}_k}$ and $\langle . , . \rangle_\mathcal{G}$. Under this RKHS framework, \cite{baker1973joint} defines the cross-covariance operator $C_k$ between $\mathcal{H}_k$ and $\mathcal{G}$ as the linear operator from $\mathcal{G}$ to $\mathcal{H}_k$ defined for all $h \in \mathcal{H}_k$ and all $g \in\mathcal{G}$ by
\begin{equation*}
\langle h , C_k g \rangle_{\mathcal{H}_k} = \cov \left( h(X_k) , g(Y) \right).
\end{equation*}
The operator $C_k$ generalizes the notion of covariance, taking into account a large spectrum of relationships between $X_k$ and $Y$ (not only linear ones). Finally, the Hilbert-Schmidt independence criterion (HSIC) is defined by \cite{gretton2005measuring} as the Hilbert-Schmidt norm of the operator $C_k$: 
\begin{equation}
\HSIC(X_k, Y)_{\mathcal{H}_k, \mathcal{G}} = \Vert C_k \Vert_{\HS}^2 = \displaystyle \sum_{i,j} \, \langle u_i ,C_k v_j \rangle_{\mathcal{H}_k}^2,  
\label{defHSIC}
\end{equation} 
where $\left( u_i \right)_{i \geq 0}$ and $\left( v_j \right)_{ j \geq 0}$ are respectively orthonormal basis of ${\mathcal{H}} _k $ and $\mathcal {G} $.

\begin{remark}
In the following, the notation $ \HSIC(X_k,Y)_ {\mathcal{H}_k, \mathcal{G}}$ is replaced by $\HSIC(X_k, Y)$ in order to lighten the expressions.
\end{remark}

\begin{sloppypar} 
Authors of \cite{gretton2005measuring} show that the HSIC measure between an input $X_k$ and the output $Y$ can be expressed using the kernels $l_k$ and $l$ in a more convenient form:
\begin{align}
\HSIC(X_k,Y) &= \mathbb{E} \left[ l_k (X_k , X'_k) l(Y , Y') \right] + \mathbb{E} \left[ l_k \left(X_k , X'_k \right) \right] \mathbb{E} \left[ l \left( Y , Y' \right) \right] \label{formHSIC} \\ &- 2 \mathbb{E} \left[  \mathbb{E} \left[ l_k \left(X_k , X'_k \right) \mid X_k \right] \mathbb{E} \left[ l \left( Y , Y' \right) \mid Y \right] \right], \nonumber 
\end{align}
where $(X'_1, \ldots ,X'_d)$ is an independent and identically distributed copy of $(X_1, \ldots ,X_d)$ and $ Y' = \mathcal{M} \left( X'_1, \ldots, X'_d \right)$. 
\end{sloppypar}

\medskip
\textbf{Independence characterization.} 
{To ensure equivalence between $\HSIC$ nullity and independence, the kernels $l_k$ and $l$ must belong to the specific class of \emph{characteristic kernels} \cite{szabo2018characteristic}. A most commonly used characteristic kernel for real variables is the Gaussian kernel, which is defined for a pair of variables} $(z, z') \in \mathbb{R}^q \times \mathbb{R}^q$ by
\begin{equation}
k_{\lambda}(z ,z') = exp \left( - \lambda \Vert z - z' \Vert_2^2 \right), 
\label{rbfkernel}
\end{equation} 
where $\lambda$ is a positive real parameter (fixed) and $\Vert . \Vert_2$ is the euclidean norm in $\mathbb{R}^q$.   

\begin{remark}
Despite that theoretically $\HSIC(X_k,Y) = 0$ is equivalent to the independence between $X_k$ and $Y$, a good choice of the kernel widths is required in practice. Indeed, a wise choice of these parameters guarantees a better behavior of HSIC estimators and better properties of the associated independence tests. Unfortunately, the best choice is unknown in practice, it depends on the joint density of $(X_k,Y)$. For this, intrinsic characteristics of these random variables are usually used. In particular, two main options are usually adopted in practice for the adjustment of $\lambda$ in Equation \eqref{rbfkernel}: whether the inverse of empirical variance of $z$, or the inverse of empirical median of $\Vert z - z' \Vert_2^2$ \cite{de2016new,sugiyama2011least,zhang2012kernel}. In the sequel, we refer to \textit{Standardized Gaussian kernel} as the one with $\lambda$ being the empirical variance. Note that, some existing works such as \cite{sugiyama2012kernel} propose methods based on cross-validation to suitably select widths. Very recently, \cite{albert2022adaptive} proposed aggregated HSIC-based tests: a well-chosen collection of HSIC tests is aggregated through a unique independence test to improve the power.   
\label{rem22}
\end{remark}

\subsubsection{Statistical estimation}
\label{Stat-estimation}

In this paragraph, we present HSIC estimators, as well as their characteristics. As a reminder, we assume that we have a \textit{n}-sample of independent realizations $\left (\mathbf{X}^{(i)}, Y^{(i)} \right)_ {1 \leq i \leq n}$ of the inputs/output couple $(\mathbf{X}, Y)$, where $\mathbf{X} = (X_1, \ldots ,X_d)$. 

\medskip
\textbf{Monte Carlo estimation.} From Equation \eqref{formHSIC}, authors of \cite{gretton2005measuring} propose to estimate each $\HSIC(X_k, Y)$ by
\begin{equation}
\widehat{\HSIC} (X_k,Y) = \displaystyle \frac{1}{n^2} \displaystyle \sum_{1 \leq i,j \leq n} (L_k)_{i,j} {(L)}_{i,j} + \displaystyle \frac{1}{n^4} \displaystyle \sum_{1 \leq i,j,q,r  \leq n} (L_k)_{i,j} {(L)}_{q,r}  - \displaystyle \frac{2}{n^3} \displaystyle \sum_{ 1 \leq i,j,r \leq n} (L_k)_{i,j} {(L)}_{j,r},
\label{esHSIC}
\end{equation}
where $L_k$ and $L$ are the matrices defined for all $i,j \in \lbrace 1, \ldots ,n \rbrace$ by $(L_k)_{i,j} = l_k ( X_k^{(i)} , X_k^{(j)} )$ and $(L)_{i,j} = l ( Y^{(i)} , Y^{(j)} )$. These V-statistic estimators \cite{mises1947asymptotic} (named after Richard Von Mises) can also be written in the following more compact form~\cite{gretton2005measuring}: 
\begin{equation}
\widehat{\HSIC} (X_k,Y) = \displaystyle \frac{1}{n^2} Tr(L_k H L H),
\label{esHSICtr}
\end{equation}
where $H$ is the matrix defined for all $i,j \in \lbrace 1, \ldots ,n \rbrace$ by $H_{i,j} = \delta_{i, j} - 1/n$, with $\delta_{i,j}$ the Kronecker symbol between $i$ and $j$ which is equal to $ 1 $ if $ i = j $ and $ 0 $ otherwise.

\medskip
\textbf{Characteristics of HSIC estimators}. Under the assumption of independence between $X_k$ and $Y$ and the assumption $l_k(x_k, x_k) = l(y, y) = 1$ (as in the case of Gaussian kernels), the estimator $\widehat{\HSIC}(X_k, Y)$ is asymptotically unbiased, its bias converges in $\mathcal{O} (1/n)$, while its variance converges to $0$ in $\mathcal{O} (1/n^2)$. Moreover, the asymptotic distribution of $n\times\widehat{\HSIC}(X_k, Y)$ is an infinite sum of independent $\chi^2$ random variables, which can be approximated by a Gamma law \cite{serfling1980approximation} with shape and scale parameters, respectively denoted $\gamma_k$ and $\beta_k$:
\begin{equation*}
\gamma_k \simeq \displaystyle \frac{e^2_k}{v_k} \quad \mbox{and} \quad \beta_k \simeq \displaystyle \frac{n.v_k}{e_k},
\end{equation*}   
where $e_k$ and $v_k$ respectively are the expectation and the variance of $\widehat{\HSIC} (X_k,Y)$, i.e. 
\begin{equation*}
e_k = \mathbb{E}\left[ \widehat{\HSIC} (X_k,Y) \right] \quad \mbox{and} \quad v_k = \Var \left( \widehat{\HSIC} (X_k,Y) \right).
\end{equation*}
The reader can refer to \cite{gretton2008kernel} and \cite{de2016new} for more details on $e_k$ and $v_k$ and their estimation. 

\subsubsection{Use for first-level GSA}
\label{GSA-HSIC}

Several methods based on HSIC measures have been developed for GSA1. In this section, we mention three possible {HSIC-based approaches for screening and ranking the inputs}: sensitivity indices \cite{da2015global}, asymptotic tests \cite{gretton2008kernel} and permutation tests~\cite{de2016new}. 

\medskip
\textbf{HSIC-based sensitivity indices.} These indices directly derived from HSIC measures, classify the input variables $X_1, \ldots, X_d$ by order of influence on the output $Y$. They are defined for all $k \in \lbrace 1, \ldots, d \rbrace$ by

\begin{equation}
\R2^2_{\HSIC,k} = \displaystyle \frac{\HSIC(X_k,Y)}{\sqrt{\HSIC(X_k,X_k) \HSIC(Y,Y)}}.
\label{RR2}
\end{equation}
  
The normalization in (\ref{RR2}) implies that $\R2^2_{\HSIC,k}$ is bounded and included in the range $[0,1]$, which makes its interpretation easier. In practice, $\R2^2_{\HSIC,k}$ can be estimated using a plug-in approach:

\begin{equation}
\widehat{\R2}^2_{\HSIC,k} = \displaystyle \frac{\widehat{\HSIC}(X_k,Y)}{\sqrt{\widehat{\HSIC}(X_k,X_k) \widehat{\HSIC}(Y,Y)}}.
\label{RR2es}
\end{equation}
  
\medskip
\textbf{Asymptotic tests.} The independence test between the input $X_k$ and the output $Y$ based on HSIC rejects the independence assumption (hypothesis denoted $\mathcal{H}_{0,k}$), when the p-value\footnote{The p-value of the test is the probability that, under $\mathcal{H}_{0, k}$, the test statistic (in this case, $n \times \widehat{\HSIC} (X_k, Y)$) is greater than or equal to the value observed on the data.} of the test based on the statistic $n \times \widehat{\HSIC} (X_k, Y)$ is less than a threshold $\alpha$ (in practice $\alpha$ is set at $5\%$ or $10\%$). Within the asymptotic framework, this p-value denoted $\Pk$ is approximated under $\mathcal{H}_{0,k}$ using the Gamma approximation (denoted~$\Gk$) of $n \times \widehat{\HSIC} (X_k, Y)$ law: 
\begin{equation}
\Pk \simeq 1 - \F_{\Gk} \left(n \times \widehat{\HSIC} (X_k,Y)_{obs} \right),
\label{testAS}
\end{equation}
where $\F_{\Gk}$ is the cumulative distribution function of $\Gk$ and $\widehat{\HSIC} (X_k,Y)_{obs}$ is the observed value of the random variable $\widehat{\HSIC} (X_k,Y)$.

\medskip
\textbf{Permutation tests.} Outside the asymptotic framework, independence tests based on permutation technique can be used. For this, the observed \textit{n}-sample is resampled $B$ independent times considering $B$ random permutations on the set $ \lbrace 1, \ldots , n \rbrace$, denoted $(\tau^{[b]})_{1 \leq b \leq B}$. These permutations are applied only to the vector $\mathbf{X}$ of inputs. We thus obtain $B$ bootstrap-samples $\left (\mathbf{X}^{\left(\tau^{[b]} (i) \right)}, Y^{(i)} \right) _ {1 \leq i \leq n}$. The HSIC measures computed on these samples are denoted $ \left( \widehat{\HSIC}^{[b]} \right)_{ 1 \leq b \leq B}$. The p-value (denoted $\pk$) of the test is then computed by
\begin{equation}
\pk = \displaystyle \frac{1}{B} \sum_{b = 1}^B  \mathds{1}_{\widehat{\HSIC}^{[b]} (X_k,Y) > \widehat{\HSIC} (X_k,Y)}.
\label{testB}
\end{equation}

{More details and demonstration of test properties are available in \cite{meynaoui2019PhD} (see Proposition 3.5). In addition, sequential algorithms have been recently proposed by \cite{elamar21}, to optimize the number of permutations $B$, while having reliable p-value estimation.}

\subsection{Estimation of HSIC with a sample generated from an alternative distribution}
\label{Estimators-alternative}

In this part, we first demonstrate that HSIC measures presented in Section \ref{Principle-definition}, can be expressed and then estimated using a sample generated from a probability distribution of inputs which is not their prior distribution. This sampling distribution will be called \enquote{alternative law} or \enquote{modified law}. The {statistical properties of these new HSIC estimators are also} presented. 
\subsubsection{Expression and estimation of HSIC measures under an alternative law}
\label{HSIC_alternative}

{We consider here $d$ continuous and independent random variables $\mathbf{\widetilde{X}} = (\widetilde{X}_1, \ldots, \widetilde{X}_d)$ whose densities (different from those of $X_1, \ldots, X_d$) are denoted $\widetilde{f}_1, \ldots, \widetilde{f}_d$. We assume that they have the same supports as $f_1, \ldots, f_d$. The associated output is denoted $\widetilde{Y} = \mathcal{M}(\mathbf{\widetilde{X}}) $. Finally, the density of~$\mathbf{\widetilde{X}}$ is designated by $ \widetilde{f}$. }

\smallskip
Changing the probability laws in HSIC expression is based on a technique commonly used in the context of importance sampling (see \textit{e.g.} {\cite{Hammersley1954}}). This technique consists in expressing an expectation $ \mathbb {E} \left[g(Z) \right]$, where $Z$ is a random variable with density $f_Z$, by using a random variable $\widetilde{Z}$ with density $f_{\widetilde {Z}}$ whose support is the same as that of $f_Z$. This gives the following expression for $\mathbb{E} \left[g(Z) \right]$:
\begin{equation}
\mathbb{E} \left[ g(Z) \right] = \displaystyle \int_{\SuppZ} g(z) \; f_Z (z) \; dz = \displaystyle \int_{\SuppZ} g(z) \; \displaystyle \frac{f_Z (z)}{f_{\widetilde{Z}} (z)} \; f_{\widetilde{Z}} (z)  \; dz = \mathbb{E}_{\widetilde{f}} \left[ g(\widetilde{Z}) \; \displaystyle \frac{f_Z (\widetilde{Z})}{f_{\widetilde{Z}} (\widetilde{Z})}  \right],
\label{Esmdf}
\end{equation}
where the notation $\mathbb{E}_{\widetilde{f}} \left[ h(Z) \right]$ designates the expectation of $h(Z)$ {for} $Z \sim \widetilde{f}$ and $\SuppZ$ denotes the support of $Z$.

\smallskip
The HSIC measures, formulated as a sum of expectations in Equation (\ref{formHSIC}), can then be expressed under the density $f_{\widetilde{Z}}$ by adapting Equation (\ref{Esmdf}) to more general forms of expectations. Hence, we obtain:
\begin{equation}
\HSIC (X_k,Y) = H^1_k + H^2_k H^3_k - 2 H^4_k, 
\label{formHSICmdf}
\end{equation}
where $(H^l_k)_{1 \leq l \leq 4}$ are the real numbers defined by
\begin{gather*}
H^1_k = \mathbb{E} \left[ l_k (\widetilde{X}_k , \widetilde{X}'_k) l(\widetilde{Y} , \widetilde{Y}') \displaystyle w (\widetilde{X}) w (\widetilde{X}')  \right]; \; H^2_k = \mathbb{E} \left[ l_k (\widetilde{X}_k , \widetilde{X}'_k ) w (\widetilde{X}) w (\widetilde{X}') \right]; \\
 H^3_k = \mathbb{E} \left[ l ( \widetilde{Y} , \widetilde{Y}')w (\widetilde{X}) w (\widetilde{X}')  \right] \; \text{and} \; H^4_k = \mathbb{E} \left[  \mathbb{E} \left[ l_k (\widetilde{X}_k , \widetilde{X}'_k ) w (\widetilde{X}')  \mid \widetilde{X}_k \right] \mathbb{E} \left[ l ( \widetilde{Y} , \widetilde{Y}' ) w (\widetilde{X}')  \mid \widetilde{Y} \right] \displaystyle w (\widetilde{X}) \right],
\end{gather*}
where $\mathbf{\widetilde{X}}'$ is an independent and identically distributed copy of $\mathbf{\widetilde{X}}$, $\widetilde{Y}' = \mathcal{M} (\mathbf{\widetilde{X}}')$ and $w = f/\widetilde{f}$.

\medskip
Formula (\ref{formHSICmdf}) shows that $\HSIC (X_k, Y)$ can then be estimated using a sample generated from $\widetilde{f}$, provided that $\widetilde{f}$ has the same support than the original density $f$. Thus, if we consider a $n$-sample of independent realizations $\left( \mathbf{\widetilde{X}}^{(i)} , \widetilde{Y}^{(i)} \right)_{ 1 \leq i \leq n }$, where $\mathbf{\widetilde{X}}$ is generated from $\widetilde{f}$ and $\widetilde{Y}^{(i)} = \mathcal{M} ( \widetilde{X}^{(i)})$, we propose the following V-statistic estimator of $\HSIC(X_k,Y)$:    
\begin{equation}
\widetilde{\HSIC} (X_k,Y) = \widetilde{H}^1_k + \widetilde{H}^2_k \widetilde{H}^3_k - 2 \widetilde{H}^4_k,
\label{esHSICmdf}
\end{equation}
where $( \widetilde{H}^l_k )_{ 1 \leq l \leq 4}$ are the V-statistics estimators of $(H^l_k)_{ 1 \leq l \leq 4}$.

\smallskip
\begin{prop} Similarly to Equation (\ref{esHSICtr}), this estimator can be rewritten as 
\begin{equation}
\widetilde{\HSIC} (X_k,Y) = \displaystyle \frac{1}{n^2} Tr \left( W \widetilde{L}_k W H_1 \widetilde{L} H_2 \right),
\label{TrHSICmdf}
\end{equation}
where $W$, $\widetilde{L}_k$, $\widetilde{L}$, $H_1$ and $H_2$ are the matrices defined by
\begin{gather*}
\widetilde{L}_k = \left( l_k ( \widetilde{X}_k^{(i)}, \widetilde{X}_k^{(j)} ) \right)_{1 \leq i,j \leq n }; \quad \widetilde{L} = \left( l ( \widetilde{Y}^{(i)}, \widetilde{Y}^{(j)} ) \right)_{1 \leq i,j \leq n }; \quad W = Diag \left( w(\widetilde{X}^{(i)}) \right)_{1 \leq i \leq n }; \\
H_1 = I_n - \displaystyle \frac{1}{n} U W; \quad H_2 = I_n - \displaystyle \frac{1}{n} W U; 
\end{gather*}
with $I_n$ is the identity matrix of size $n$ and $U$ the matrix filled with $1$.
\end{prop}

\smallskip
The proof of this proposition is detailed in Appendix \ref{sec:AnnexA}. 
Similarly, the sensitivity index $\R2^2_{\HSIC,k}$ can also be estimated using the sample $\left(\mathbf{\widetilde{X}}^{(i)},\widetilde{Y}^{(i)}\right)_{1 \leq i \leq n}$ by
\begin{equation}
\widetilde{\R2}^2_{\HSIC,k} = \displaystyle \frac{\widetilde{\HSIC} (X_k,Y)}{\sqrt{\widetilde{\HSIC} (X_k,X_k) \widetilde{\HSIC} (Y,Y)}}.
\label{R2mdf}
\end{equation} 

\smallskip
We can demonstrate that these new estimators have statistical properties similar to those of classical estimators. More precisely,  $\widetilde{\HSIC} (X_k, Y)$ is asymptotically unbiased and its bias converges in $\mathcal{O} (1/n)$. 
\begin{prop} Under the hypothesis of independence of $X_k$ and $Y$, and the assumption $l_k(x_k, x_k) = l(y, y) = 1$, the bias and variance of $\widetilde{\HSIC} (X_k, Y)$ are respectively:
\begin{align}
\mathbb{E} \left[ \widetilde{\HSIC} (X_k,Y) \right] - \HSIC \left( X_k,Y \right) &=  \frac{2}{n}( E_{\omega}^k -  E_{x_k,\omega})(E_{\omega}^{-k} - E_{y,\omega} ) - \frac{1}{n} (E_{\omega} - E_{x_k})(E_{\omega} - E_y) \label{BiaisHSICmdf} \\ &+  \frac{1}{n} E_{\omega}(E_{\omega} - 1) + \mathcal{O}(1/n^2), \nonumber \\ 
\Var \left[ \widetilde{\HSIC} (X_k,Y) \right] &= \frac{72(n-4)(n-5)}{n(n-1)(n-2)(n-3)} \; \mathbb{E}_{1,2} \left[ \mathbb{E}_{3,4} [ \widetilde{h}_{1,2,3,4} ]^2 \right] + \mathcal{O} (1/n^3),
\label{VarHSICmdf}
\end{align}
where
\begin{alignat*}{4}
&E_{\omega} &{}={}& \mathbb{E} \left[ \omega^2(\widetilde{X})\right], &{} &E_{x_k} &{}={}& \mathbb{E} \left[ l_k ( \widetilde{X}_k , \widetilde{X}'_k ) \omega_k(\widetilde{X}_k) \omega_k(\widetilde{X}'_k)  \right], \\
&E_y  &{}={}& \mathbb{E} \left[ l( \widetilde{Y} , \widetilde{Y}' ) \omega_{-k}(\widetilde{X}_{-k}) \omega_{-k}(\widetilde{X}'_{-k})  \right], \hspace{0.5cm} &{}  &E_{x_k,\omega} &{}={}& \mathbb{E} \left[ l_k ( \widetilde{X}_k , \widetilde{X}'_k ) \omega^2_k(\widetilde{X}_k)  \omega_k(\widetilde{X}'_k)  \right], \\
&E_{y,\omega} &{}={}& \mathbb{E} \left[ l ( \widetilde{Y} , \widetilde{Y}' ) \omega^2_{-k}(\widetilde{X}_{-k}) \omega_{-k}(\widetilde{X}'_{-k})  \right], &{} &E_{\omega}^k &{}={}& \mathbb{E} \left[ \omega_k^2(\widetilde{X}_k)\right], \\
&E_{\omega}^{-k} &{}={}& \mathbb{E} \left[ \omega_{-k}^2(\widetilde{X}_{-k})\right], &{} & \widetilde{h}_{1,2,3,4} &{}={}&  \frac{1}{4!} \sum_{(t,u,v,s)}^{(1,2,3,4)} \left[ (\widetilde{l}_k)_{t,u}  \widetilde{l}_{t,u} + (\widetilde{l}_k)_{t,u} \widetilde{l}_{v,s}  - 2 (\widetilde{l}_k)_{t,u} \widetilde{l}_{t,v} \right],        
\end{alignat*}
$\omega$, $\omega_k$ and $\omega_{-k}$ respectively denote the functions $f/\widetilde{f}$ , $f_k/\widetilde{f}_k$ and $\omega_{-k} : x_{-k} \mapsto \omega(x_1, \ldots, x_d)/\omega_k (x_k)$, with $x_{-k}$ being the vector extracted from $(x_1,\ldots ,x_d)$ by removing the~k-th~coordinate. Moreover, $\widetilde{X}'_{-k}$ is an independent and identically distributed copy of $\widetilde{X}_{-k}$ and $(\widetilde{l}_k)_{p,q}$, $\widetilde{l}_{p,q}$ respectively denote $l_k ( \widetilde{X}^{(p)}_k , \widetilde{X}^{(q)}_k )$, $l ( \widetilde{Y}^{(p)} , \widetilde{Y}^{(q)})$. Finally, $\sum_{(t,u,v,s)}^{(1,2,3,4)}$ is the sum over all permutations $(t,u,v,s)$ of $(1,2,3,4)$ and $\mathbb{E}_{p,q}$ is the expectation only with respect to $\mathbf{X}_p$ and $\mathbf{X}_q$.
\end{prop}

{One can also prove that the distribution of $n \times \widetilde{\HSIC} (X_k, Y) $ can be approximated by a Gamma law, whose parameters $\widetilde{\gamma}_k$ and $\widetilde{\beta}_k$ are given by $\widetilde{\gamma}_k = \varepsilon^2_k/\vartheta_k$ and $ \widetilde{\beta}_k = n \vartheta_k/\varepsilon_k$, where $\varepsilon_k$ and $\vartheta_k $ are respectively the expectation and variance of $ \widetilde{\HSIC} (X_k, Y)$.
Proofs of all the propositions are provided in \cite{meynaoui2019PhD}, as well as unbiaised estimators of bias and variance.}

\begin{remark}
From a practical point of view, the greater $(\Var ( \omega_k (\widetilde{X}_k) )_{1 \leq k \leq d}$, the greater the number of simulations required to accurately estimate $ (\HSIC (X_k,Y))_{1 \leq k \leq d}$. It is therefore highly recommended to check that $(\Var ( \omega_k (\widetilde{X}_k) )_{1 \leq k \leq d}$ are finite. For instance, in the case of densities with compact supports, it is enough to check that $(\omega_k)_{1 \leq k \leq d}$ are finite on their supports. 
\end{remark}

\subsubsection{Illustration on an analytical example}
\label{example_num}

{To illustrate the statistical properties of $\widetilde{\HSIC}$, we consider a numerical application inspired from Ishigami's model \cite{ishigami1990importance} and defined on $[0,1]^3$ by}
\begin{equation}
\mathcal{M} (X_1,X_2,X_3) = \sin(X_1) + 1.5 \sin^2 (X_2) + 0.5 X_3^4 \sin(X_1),
\label{Ishimdf}
\end{equation}
where the inputs $X_1$, $X_2$ and $X_3$ are assumed to be independent and follow a triangular distribution with a mode equal to $0.5$.

\begin{figure}[h!]
\includegraphics[scale=0.15]{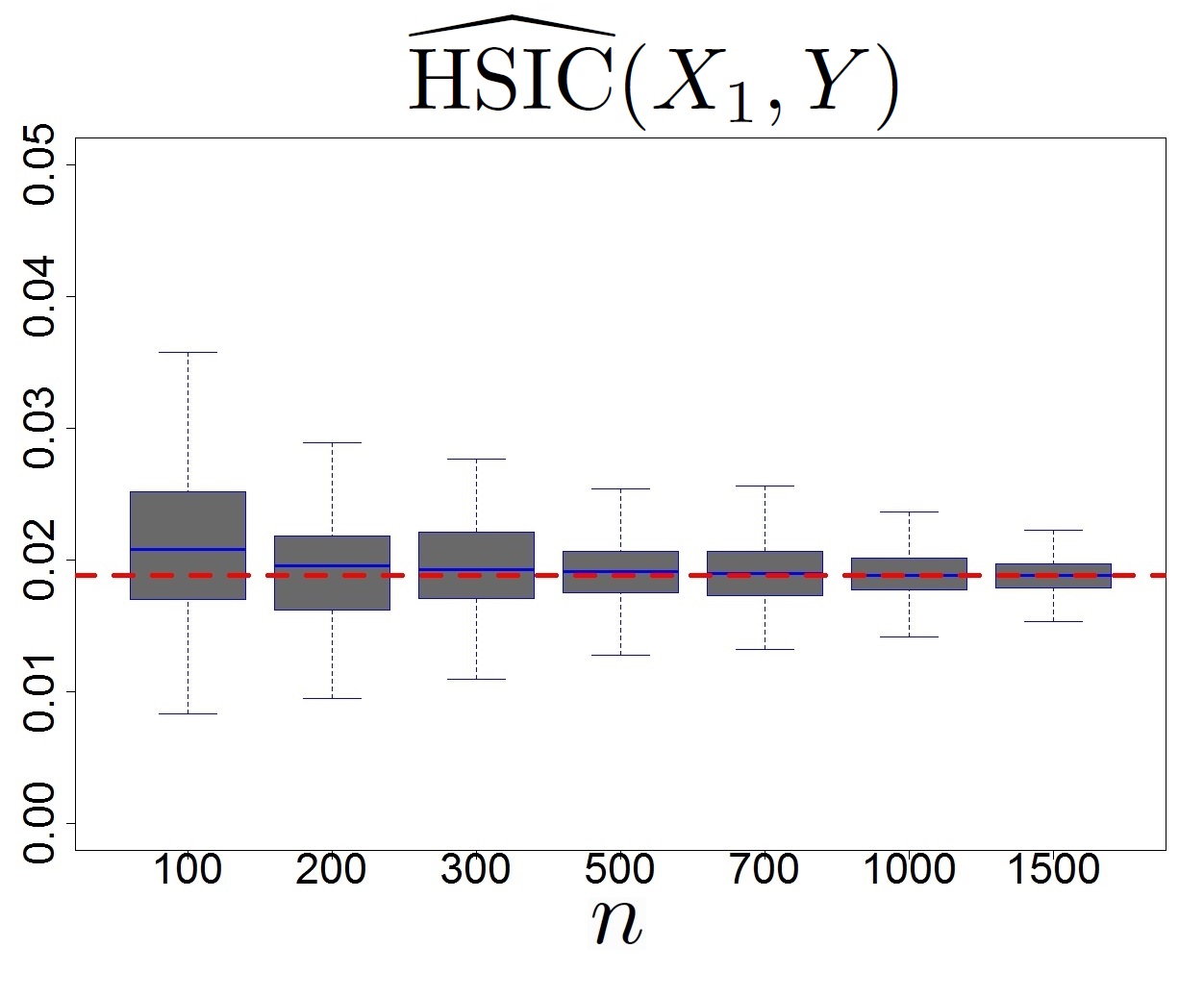} 
\includegraphics[scale=0.15]{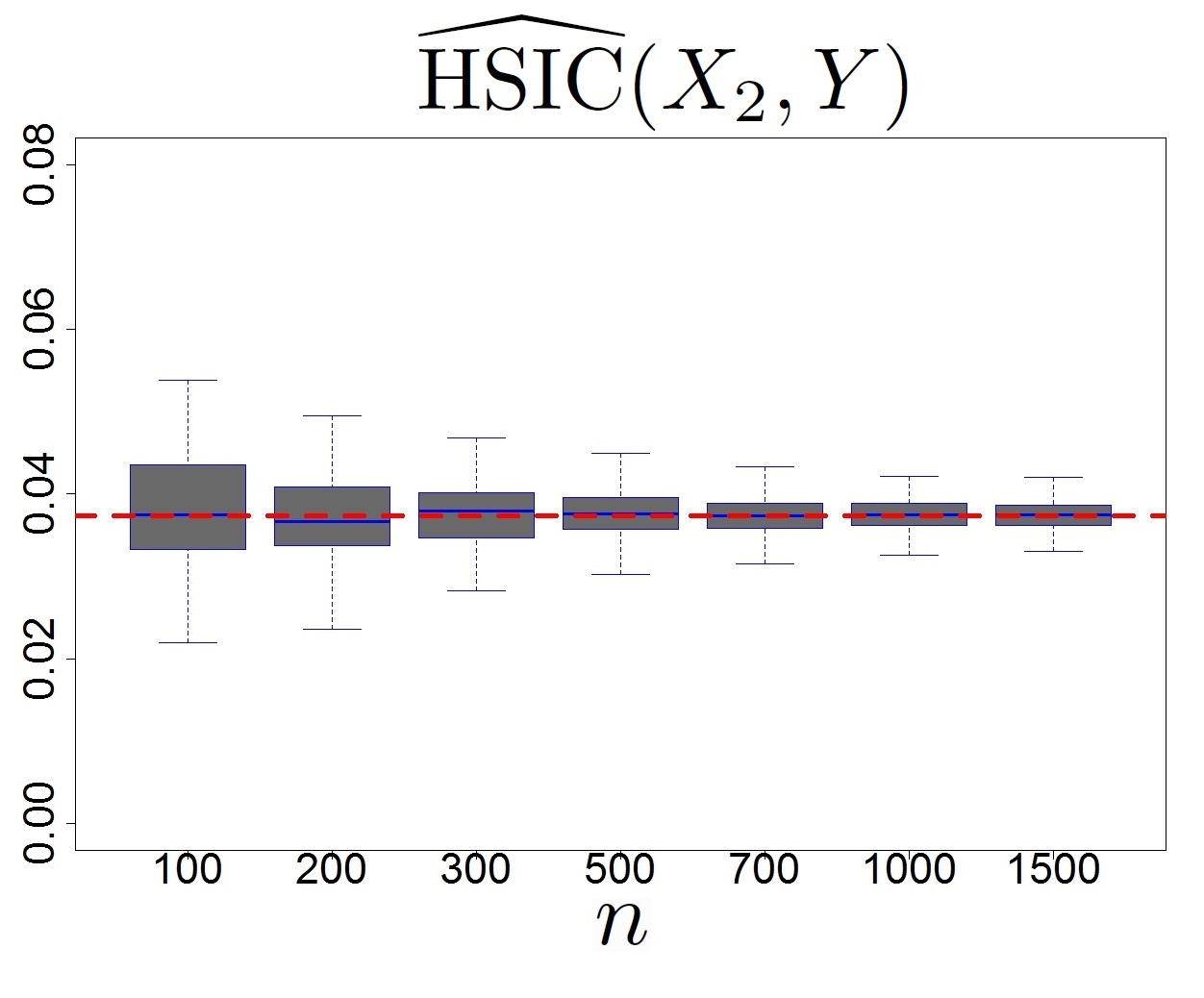}
\includegraphics[scale=0.15]{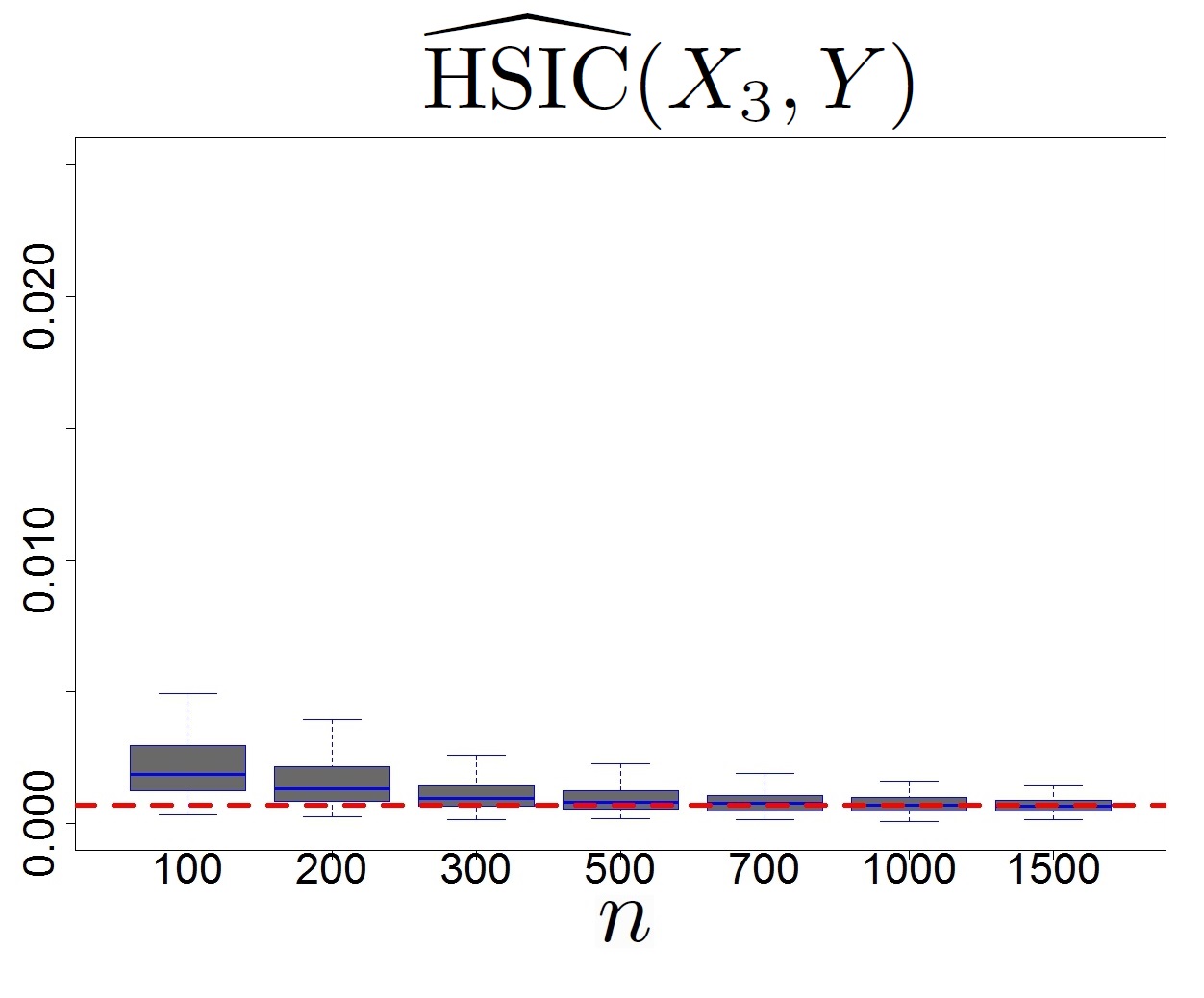}
\centering
\includegraphics[scale=0.15]{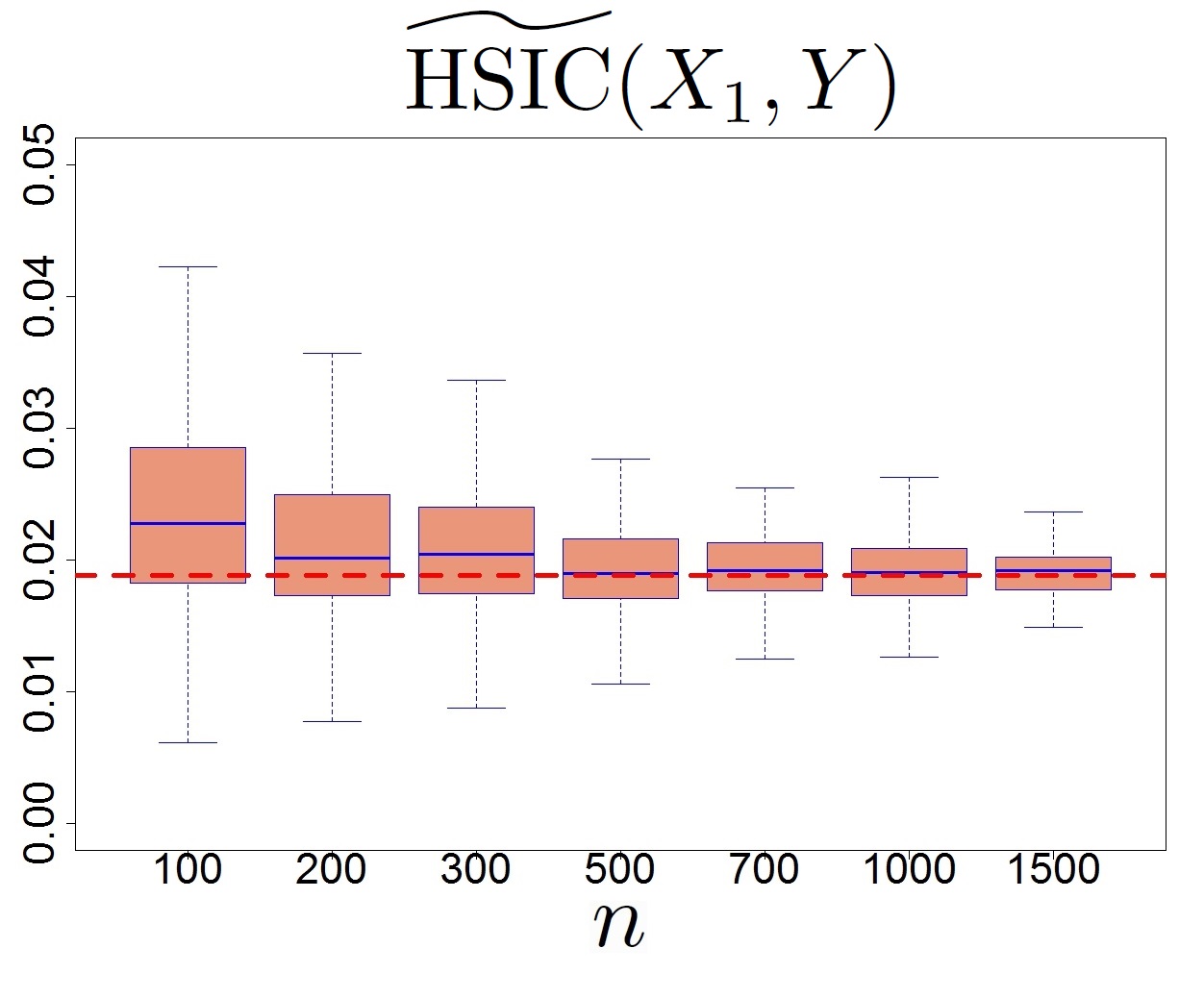} 
\includegraphics[scale=0.15]{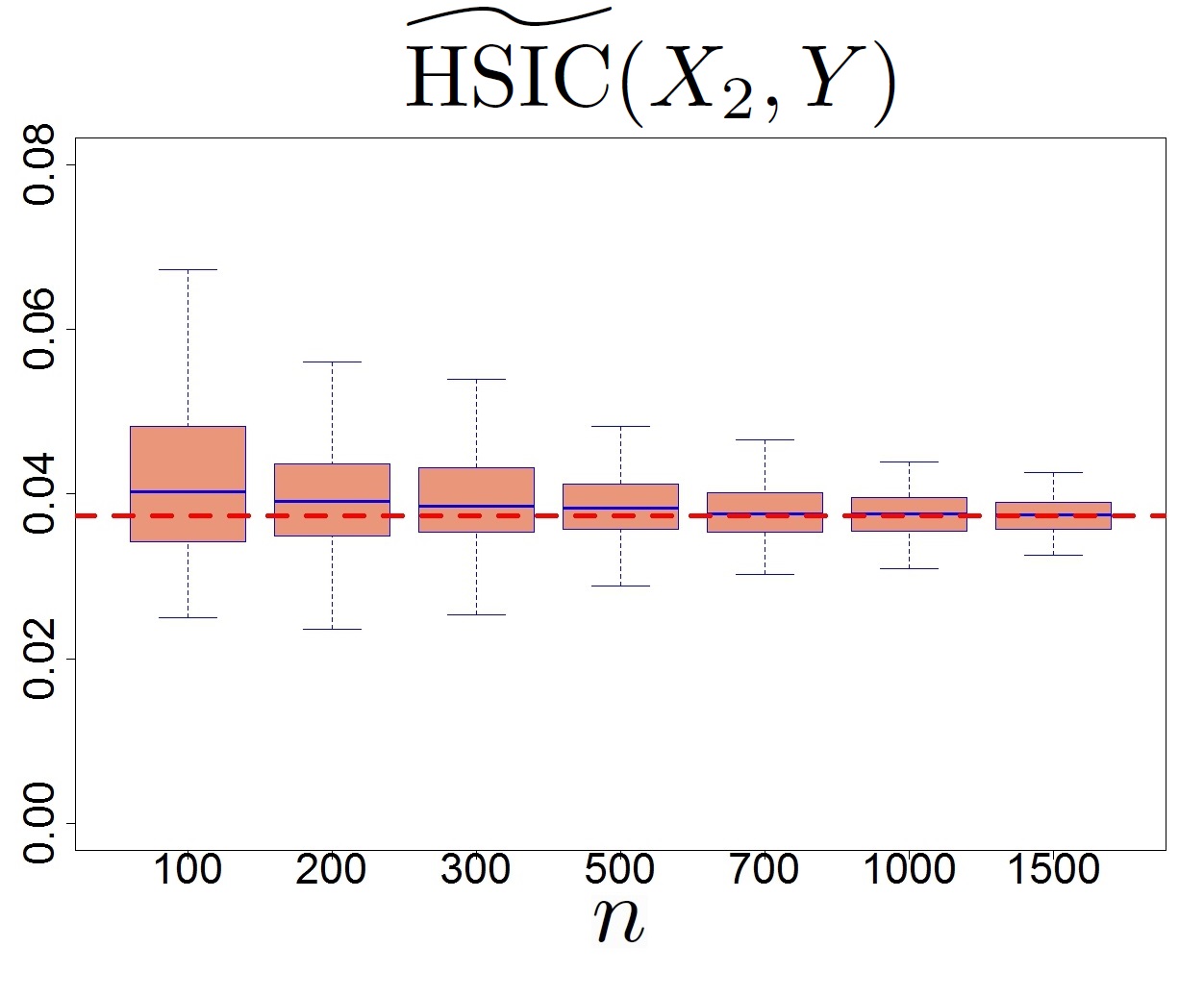}
\includegraphics[scale=0.15]{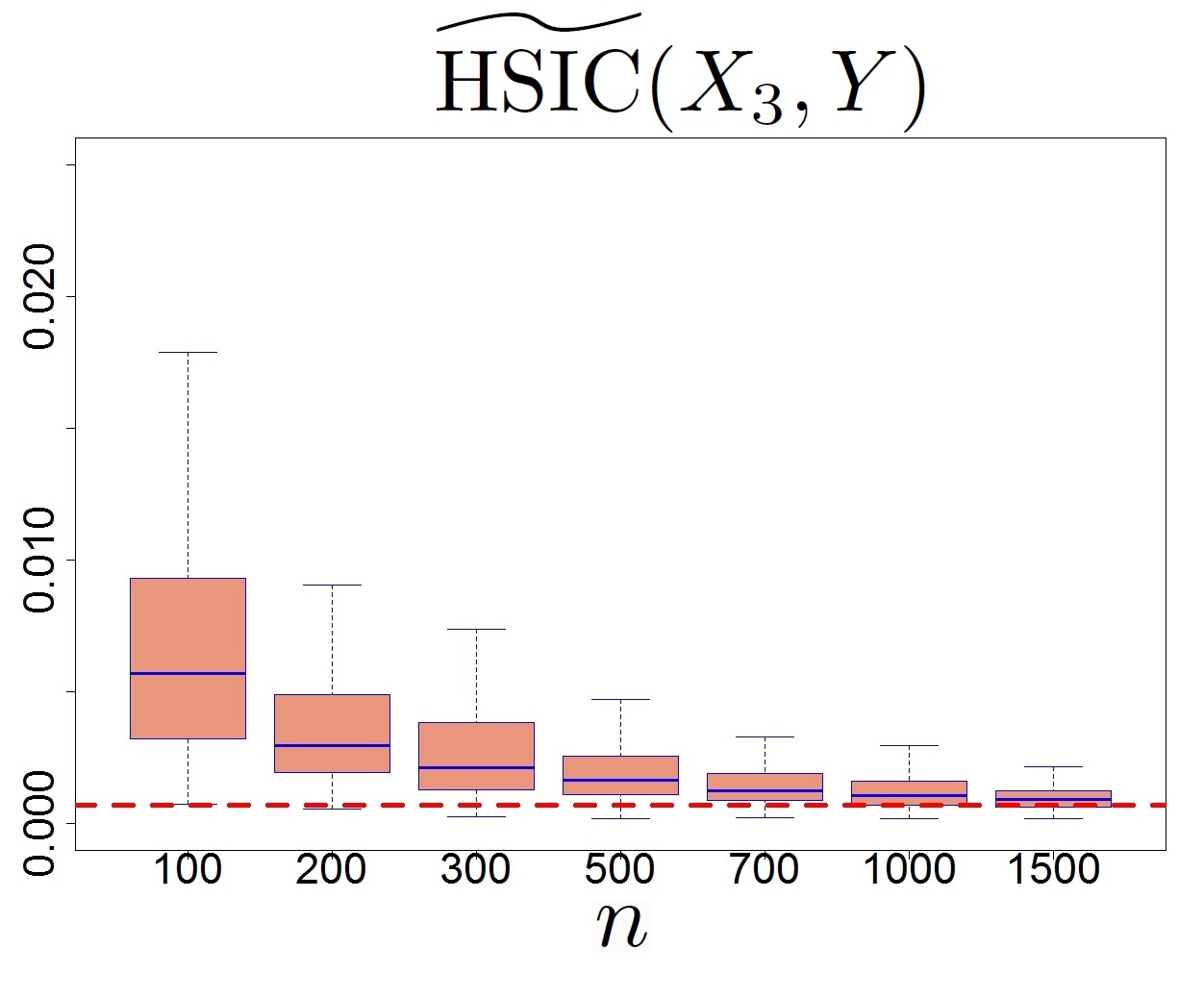}
\caption{{Model $\mathcal{M}$ -- Convergence plots of the estimators $\widehat{\HSIC} (X_k,Y)$ and $\widetilde{\HSIC} (X_k,Y)$, according to the sample size $n$. Theoretical values are represented in red dashed lines.}}
\label{VvstatHSIC}
\end{figure}

\smallskip
We consider HSIC measures based on Standardized Gaussian kernel (see Remark \ref {rem22}). {To estimate them, we suppose that we have Monte Carlo samples of independent inputs generated from a uniform distribution on $[0,1]^3$ (modified law). For each sample of size $ n = 100 $ to $ n = 1500 $, the estimation process is repeated $200$ times, with independent random samples.}
The convergence of $\widetilde{\HSIC}(X_k,Y)$ estimators is illustrated by Figure \ref{VvstatHSIC}. Results for $\widehat{\HSIC}(X_k,Y)$ computed with samples generated from the original law (namely triangular) are also given and theoretical values are represented in red dotted lines. We observe that for small sample sizes ($n <500$), modified estimators $\widetilde{\HSIC} (X_k,Y)$ have more bias and variance than $\widehat{\HSIC}(X_k, Y)$ estimators. But, from size $n = 700$, both estimators have similar behaviors. 

\smallskip
In addition, to assess the convergence of ranking, the sensitivity indices $ R^2_{\HSIC,k}$ are estimated from $\widetilde{\HSIC}(X_k,Y)$ with Equation (\ref{R2mdf}). The inputs are ranked by decreasing indices and the resulting correct ranking rates are given by Table \ref{tauxbon}. Even for small sample sizes (\textit{e.g} $ n = 200 $), the modified estimators $\widetilde{\R2}^2_{\HSIC}$ have good ranking ability.

\newcolumntype{D}[1]{>{\centering}p{#1}}
\begin{table}[!h]
\begin{center}
\begin{tabular}{|D{1.8cm}|D{1.8cm}|D{1.8cm}|D{1.8cm}|}\hline
$n = 100$ & $n = 200$ & $n = 300$ & $n \geqslant 500$ \tabularnewline\hline \hline
$88 \%$ & $93.5\%$ & $97\%$& $100\%$ \tabularnewline\hline
\end{tabular} 
\caption{{Model $\mathcal{M}$ -- Good ranking rates of inputs based on $\widetilde{\R2}^2_{\HSIC}$, for different sample sizes $n$.}}
\label{tauxbon}
\end{center}
\end{table}

\section{New methodology for second-level GSA}
\label{NEW-GSA2}

We consider that the input probability distributions, $\mathbb{P}_{X_1}, \ldots , \mathbb{P}_{X_d}$, are uncertain. These uncertainties on $\mathbb{P}_{X_1}, \ldots , \mathbb{P}_{X_d}$ are modeled by probability distributions, respectively denoted $\mathbb{P}_{\mathbb{P}_{X_1}}, \ldots, \mathbb{P}_{\mathbb{P}_{X_d}}$. We also assume that the distributions $\mathbb{P}_{\mathbb{P}_{X_1}}, \ldots, \mathbb{P}_{\mathbb{P}_{X_d}}$ are independent and that all possible input distributions have a common support, which the set of all possible input values.
Each assumed joint distribution $\mathbb{P}_{\mathbf{X}} = \mathbb{P}_{X_1} \times \ldots \times \mathbb{P}_{X_d}$ of inputs yields potentially different results of 1$^{\text{st}}$-level global sensitivity analysis (GSA1). This impact must be quantified by GSA2. Based on GSA2 results, the probability distributions of inputs could be separated into two groups: those which significantly modify GSA1 results and those whose influence is negligible. Subsequently, the probability distributions with a small impact can be set to a reference distribution and the efforts of characterization will be focused on the most influential ones to improve their knowledge (strategy of uncertainty reduction).

\subsection{Issues raised by GSA2}
\label{GSA2-Issues}

We present in the following the main steps for GSA2 realization and some related technical issues. Our approach is based on the extension of HSIC measures for non-vectorial data. The idea is to define 2$^{\text{nd}}$-level sensitivity indices between input distributions $\mathbb{P}_{X_1}, \ldots, \mathbb{P}_{X_d}$ and GSA1 results. To do so, we first characterize GSA1 results. This means that we associate to each possible input distribution $\mathbb{P}_{\mathbf{X}} = \mathbb{P}_{X_1} \times \ldots \times \mathbb{P}_{X_d}$, a mathematical quantity denoted $\Res$ representing the associated GSA1 results. To choose this quantity of interest, we propose the following options, all based on HSIC (see Section \ref{GSA-HSIC}):

\begin{itemize}
\item \textbf{vector of sensitivity indices} $\mathbf{\R2^2_{\HSIC}} =(\R2^2_{\HSIC,1}, \ldots , \R2^2_{\HSIC,d})$; 

\item \textbf{ranking of inputs} $X_1 , \ldots , X_d$ using the indices $\R2^2_{\HSIC,1}, \ldots , \R2^2_{\HSIC,d}$ \textbf{.} This quantity of interest $\Res$ is a permutation on the set $\lbrace 1, \ldots ,d \rbrace$,  verifying that $\Res (k) = j$ if and only if the variable $X_j$ is the $k$-th in the ranking; 

\item \textbf{vector of p-values from asymptotic independence tests} $\mathbf{\Pval} = (\P1, \ldots , \Pd)$;

\item \textbf{vector of p-values from permuted tests} $\mathbf{\pval} = (\p1, \ldots , \pd)$. 
\end{itemize}

Thanks to the kernel trick, we build 2$^{\text{nd}}$-level HSIC measures between the probability distributions $\mathbb{P}_{X_1}, \ldots ,\mathbb{P}_{X_d}$ and the quantity of interest $\Res$. Assume that $l_{\D1}, \ldots ,l_{\Dd}$ and $l_{\Res}$ are RKHS kernels respectively associated to $\mathbb{P}_{X_1}, \ldots ,\mathbb{P}_{X_d}$ and $\Res$. Some examples of these kernels are provided in Section \ref{RKHS-choices}. We define similarly to Equation \eqref{formHSIC}, the 2$^{\text{nd}}$-level HSIC measure between $\mathbb{P}_{X_k}$ and $\Res$ as
\begin{align}
\HSIC(\mathbb{P}_{X_k},\Res) &= \mathbb{E} \left[ l_{\Dk} (\mathbb{P}_{X_k}, \mathbb{P}'_{X_k}) l_{\Res}(\Res , \Res') \right] + \mathbb{E} \left[ l_{\Dk} ( \mathbb{P}_{X_k} , \mathbb{P}'_{X_k} ) \right] \mathbb{E} \left[ l_{\Res} \left( \Res , \Res' \right) \right] \nonumber \\ &- 2 \mathbb{E} \left[  \mathbb{E} \left[ l_{\Dk} (\mathbb{P}_{X_k} , \mathbb{P}'_{X_k} ) \mid \mathbb{P}_{X_k} \right] \mathbb{E} \left[ l_{\Res} \left( \Res , \Res' \right) \mid \Res \right] \right],
\label{formHSIC2}
\end{align} 
where $(\mathbb{P}'_{X_1}, \ldots, \mathbb{P}'_{X_d} )$ is an independent and identically distributed copy of $(\mathbb{P}_{X_1}, \ldots, \mathbb{P}_{X_d})$ and $\Res'$ the GSA1 results associated to $(\mathbb{P}'_{X_1}, \ldots, \mathbb{P}'_{X_d})$. The GSA2 indice between $\mathbb{P}_{X_k}$ and $\Res'$ is then defined as
\begin{equation}
\R2^2_{\HSIC} (\mathbb{P}_{X_k} , \Res) = \displaystyle \frac{\HSIC(\mathbb{P}_{X_k},\Res) }{\sqrt{\HSIC(\mathbb{P}_{X_k},\mathbb{P}_{X_k})\HSIC(\Res,\Res)}}.
\label{formR2HSIC2}
\end{equation}

The estimation of $\R2^2_{\HSIC} (\mathbb{P}_{X_k} , \Res)$  requires a $n_1$-sample $(\mathbb{P}_{\mathbf{X}}^{(i)},\mathcal{R}^{(i)})_{ 1 \leq i \leq n_1}$ of $(\mathbb{P}_{\mathbf{X}},\mathcal{R})$. However, the quantities of interest $\Res^{(i)}$ are not directly observable, they need to be estimated. To do so, a straightforward double-loop approach could be considered. The outer loop entails to generate the $n_1$-sized sample of input distribution. On the flip side, the inner loop involves two steps. A $n_2$-sized sample $(X_1^{(i,j)}, \ldots ,X_d^{(i,j)} )_{1 \leq j \leq n_2}$ is first generated according to each distribution $\mathbb{P}_{\mathbf{X}}^{(i)}$, before computing the corresponding outputs $(Y^{(i,j)})_{1 \leq j \leq n_2}$. This allows to estimate the quantity of interest $\mathcal{R}^{(i)}$ associated to each input distribution~$\mathbb{P}_{\mathbf{X}}^{(i)}$. At the end, 2$\textsuperscript{nd}$-level HSIC can be estimated by
\begin{equation}
\widehat{\HSIC} (\mathbb{P}_{X_k},\mathcal{R}) = \displaystyle \frac{1}{n_1^2} \Tr(L_{\Dk} H L_{\mathcal{R}} H),
\label{esHSICtr(2)}
\end{equation}
where $L_{\Dk}$ and $L_{\mathcal{R}}$ are the matrices defined for all $(i,j)$ in $\lbrace 1, \ldots ,n_1 \rbrace^2$ as $( L_{\Dk} )_{i,j} = l_{\Dk} ( \mathbb{P}^{(i)}_{X_k} , \mathbb{P}^{(j)}_{X_k} )$ and $(L_{\mathcal{R}})_{i,j} = l_{\mathcal{R}} ( \mathcal{R}^{(i)} , \mathcal{R}^{(j)} )$. In addition, the matrix $H$ is defined as in Equation (\ref{esHSICtr}). Finally, according to \eqref{formR2HSIC2}, 2$^{\text{nd}}$-level $\R2^2_{\HSIC}$ indices can be estimated by 
\begin{equation}
\widehat{\R2}^2_{\HSIC} (\mathbb{P}_{X_k},\mathcal{R}) = \displaystyle \frac{\widehat{\HSIC}(\mathbb{P}_{X_k},\mathcal{R})}{\sqrt{\widehat{\HSIC}(\mathbb{P}_{X_k},\mathbb{P}_{X_k}) \widehat{\HSIC}(\mathcal{R},\mathcal{R})}}.
\label{RR2(2)}
\end{equation}

Consequently, the Monte Carlo double-loop approach requires a total of $n_1 n_2$ simulations. This approach is therefore not tractable for CPU-time expensive simulators, even for reasonable sample sizes $n_1$ and $n_2$. To overcome this issue and reduce the number of simulator-calls, we propose in the following a single-loop approach only requiring $n_2$ simulations, and allowing to consider a larger sample of input distribution.

\subsection{Algorithm for computing 2$^{\text{nd}}$-level sensitivity indices with a Monte Carlo single-loop}
\label{1-loop-algorithm}
We provide here a single-loop algorithm for estimating the 2$^{\text{nd}}$-level HSIC measures and indices, respectively defined in Equations \eqref{esHSICtr(2)} and \eqref{RR2(2)}. To do so, the inputs are generated according to a unique and known probability distribution, denoted $\bar{\mathbb{P}}_{\mathbf{X}} = \bar{\mathbb{P}}_{X_1} \times \ldots \times \bar{\mathbb{P}}_{X_d}$. We assume that this distribution has a density $\bar{f} : (x_1 , \ldots , x_d) \mapsto \bar{f}_1 (x_1) \times \ldots \times \bar{f}_d (x_d)$, and that all possible input distributions also have densities. The procedure is detailed in Algorithm \ref{algo_GSA2}.

\begin{algorithm}[h] 
\caption{\textit{GSA2 with a Monte Carlo single-loop}}\label{algo_GSA2}
\flushleft\textbf{Input:} The probability density $\bar{f}$ and an observed $n_2$-sized sample $\bar{\mathbf{X}} = \left( \mathbf{X}^{(1)}, \ldots,  \mathbf{X}^{(n_2)} \right)$.
\justify
\begin{enumerate} 
\item \label{Unique-Sample} \textbf{Build a unique $n_2$-sized sample $\mathcal{E}$ of inputs/output.}

We compute the output sample $\bar{Y} = \left( Y^{(1)}, \ldots, Y^{(n_2)} \right)$ associated to $\bar{\mathbf{X}}$. The inputs/output sample is denoted $\mathcal{E} = \left( \bar{\mathbf{X}}, \bar{Y} \right)$.

\item \label{Multiple-GSA1} \textbf{Perform $n_1$ GSA1 using only $\mathcal{E}$.}

We draw a $n_1$-sized sample $ \mathbb{P}_{\mathbf{X}}^{(1)}, \ldots, \mathbb{P}_{\mathbf{X}}^{(n_1)}$ of input distributions.
 Then, we estimate all the GSA1 results $\mathcal{R}^{(i)}$ associated to each distribution $\mathbb{P}_{\mathbf{X}}^{(i)}$ \underline{using only} $\mathcal{E}$. The options proposed for $\mathcal{R}^{(i)}$ in Section \ref{GSA2-Issues} are distinguished: 

\begin{enumerate}
\item \label{QI-R2} \textbf{Estimate the vector} $\mathcal{R}^{(i)} =(\R2^{2,(i)}_{\HSIC,1}, \ldots , \R2^{2,(i)}_{\HSIC,d})$ \textbf{of sensitivity indices.} 

The vector coordinates are estimated using Equation \eqref{R2mdf}, with the alternative sample $\mathcal{E} = \left( \bar{\mathbf{X}}, \bar{Y} \right)$.

\item \label{QI-RK}  \textbf{Rank the inputs} $X_1 , \ldots , X_d$ \textbf{using the indices} $\R2^2_{\HSIC,1}, \ldots , \R2^2_{\HSIC,d}$ \textbf{.} \newline
The ranking is obtained by ordering the coordinates of the vectors estimated in Option \ref{QI-R2}.
 
\item \label{Pval-asy} \textbf{Estimate the vector} $\mathcal{R}^{(i)} = \left(\Pval_1^{(i)}, \ldots , \Pval_d^{(i)}\right)$ \textbf{of p-values associated with asymptotic independence tests.} \newline
Each $\Pval_k^{(i)}$ is estimated using the properties of the modified estimators:      
\begin{equation}
\widetilde{\Pval}_k^{(i)} \simeq 1 - \widetilde{\F}_{\Gk} \left( n_2 \times \widetilde{\HSIC} (X_k^{(i)},Y)_{obs} \right),
\label{testASmdf}
\end{equation}
where $\tilde{\F}_{\Gk}$ denotes the Gamma distribution approximating the asymptotic distribution of $n_2 \times \widetilde{\HSIC} (X_k^{(i)},Y)$.

\item \label{Pval-perm} \textbf{Estimate the vector} $\mathcal{R}^{(i)} = \left(\pval_1^{(i)}, \ldots , \pval_d^{(i)}\right)$ \textbf{of p-values associated with permutation independence tests.} \newline 
Keeping the same notations of Equation \eqref{testB}, each $\pval_k^{(i)}$ is estimated as
\begin{equation}
\widetilde{\pval}_k^{(i)} = \displaystyle \frac{1}{B} \sum_{b = 1}^B  \mathds{1}_{\widetilde{\HSIC}^{[b]} (X_k^{(i)},Y) > \widetilde{\HSIC} (X_k^{(i)},Y)}.
\label{testBmdf}
\end{equation}
\end{enumerate}

\item\label{Estim-ind-GSA2} \textbf{Compute 2$^{\text{nd}}$-level sensitivity indices.}

Each indice $\R2^2_{\HSIC} \left(\mathbb{P}_{X_k},\Res \right)$ is estimated from the sample $\left(\mathbb{P}^{(i)}_{\mathbf{X}},\widetilde{\mathcal{R}}^{(i)}\right)_{ 1 \leq i \leq n_1}$ and using Equations \eqref{esHSICtr(2)} and \eqref{RR2(2)}. \newline
\end{enumerate}
\label{AlgoGSA2}
\end{algorithm}

\subsection{Choice of characteristic kernels for probability distributions and for quantities of interest}
\label{RKHS-choices}

Step \ref{Estim-ind-GSA2} of Algorithm \ref{AlgoGSA2} involves a choice of kernel, according to the quantities of interest $\mathcal{R}$. A kernel for probability distributions is also required. Some examples of suitable characteristic RKHS kernels are provided in the following.\\

\textbf{Characteristic RKHS kernel for probability distributions.} 
The definition of RKHS kernels between distributions is based on the Maximum Mean Discrepancy (MMD), introduced in \cite{gretton2012kernel}. Let $\mathbb{Q}_1$ and $\mathbb{Q}_2$ be two distributions with a common support and $K$ be a RKHS kernel on this support. The distance MMD between $\mathbb{Q}_1$ and $\mathbb{Q}_2$ is defined as
\begin{equation}
\text{MMD}_K (\mathbb{Q}_1 , \mathbb{Q}_2) = \mathbb{E}\left[ K(Z_1,Z'_1) \right] - 2 \mathbb{E} \left[K(Z_1,Z_2) \right] + \mathbb{E} \left[ K(Z_2,Z'_2) \right], 
\label{MMD}
\end{equation}
where $Z_1$ and $Z'_1$ (respectively $Z_2$ and $Z'_2$) are independent random variables with common distribution $\mathbb{Q}_1$ (respectively $\mathbb{Q}_2$). From the MMD distance, we consider the radial RKHS distribution kernel defined as
\begin{equation}
 l_{\mathcal{D}} \left(\mathbb{Q}_1 , \mathbb{Q}_2\right) = \exp \left( - \lambda \MMD^2_{K} \left(\mathbb{Q}_1 , \mathbb{Q}_2\right) \right), 
\label{kernelMMD}
\end{equation}
where $\lambda$ is a positive real parameter (fixed). The latter RKHS kernel is characteristic, regardless of $\lambda$ value. This property results from \citep[Theorem 4]{sriperumbudur2009kernel}. However, the parameter $\lambda$ needs to be properly calibrated for well behaved estimators. Back to our case, we define the kernel $l_{\mathcal{D}_k}$ in Equation \eqref{esHSICtr(2)} by choosing $K = l_k$ and $\lambda = 1/s_k^2$ where
$$s^2_k = \frac{1}{n_1^2} \sum_{i = 1}^{n_1} \MMD_{l_k}^2 \left( \mathbb{P}^{(i)}_{X_k}, \mathbb{P}^{\dagger}_{X_k} \right),$$ 
with $\mathbb{P}^{\dagger}_{X_k} = 1/n_1 \sum_{i=1}^{n_1} \mathbb{P}^{(i)}_{X_k}$ is the uniformly weighted mixture distribution of $\mathbb{P}^{(1)}_{X_k}, \ldots, \mathbb{P}^{(n_1)}_{X_k}$.\\

\textbf{Characteristic RKHS kernel for permutations as quantity of interest.} 
When considering Option \ref{QI-RK} in Algorithm \ref{AlgoGSA2}, we propose the use of the Mallows kernel, brought to light by \cite{jiao2018kendall}. The Mallows kernel is shown to  be universal (and characteristic) by \cite{mania2018kernel}. To define this kernel, we first introduce the \textit{number of discordant pairs} between two $\{ 1, \ldots, d \}$-permutations $\sigma_1$ and $\sigma_2$ as
\begin{equation}
n_D (\sigma_1 , \sigma_2) = \sum_{1 \leq r < s \leq d} \left[ \mathds{1}_{\lbrace \sigma_1 (r) < \sigma_1 (s) \rbrace} \mathds{1}_{\lbrace \sigma_2(r) > \sigma_2(s) \rbrace} + \mathds{1}_{\lbrace \sigma_1(r) > \sigma_1(s) \rbrace} \mathds{1}_{\lbrace \sigma_2(r) < \sigma_2(s) \rbrace} \right]. 
\end{equation}
The Mallows kernel $K_M$ between $\sigma_1$ and $\sigma_2$ is then defined as
\begin{equation}
K_M (\sigma_1 , \sigma_2) = \exp \left( - \lambda n_D ( \sigma_1 , \sigma_2) \right),
\label{MallowsKernel}
\end{equation} 
where $\lambda$ is a positive real. From a numerical standpoint and based on a sample $\sigma^{(1)}, \ldots, \sigma^{(n_1)}$, we suggest to take $\lambda$ as the inverse of the arithmetic mean of $\left\{ n_D (\sigma^{(i)} , \sigma^{(j)}) \mbox{ with } 1 \leq i < j \leq n_1 \right\}$.\\

\textbf{Characteristic RKHS kernel for real vectors as quantities of interest.} When either Option \ref{QI-R2}, \ref{Pval-asy} or \ref{Pval-perm} is selected, we can simply use the Standardized Gaussian kernel (see Remark \ref{rem22}).

\subsection{Possibilities for the unique sampling distribution}
\label{unique-law}

This section deals with the choice of the drawing density $\bar{f}$ in Algorithm \ref{AlgoGSA2}. Since the inputs are independent, it boils down to choosing each martingale density $\bar{f}_k$. We recall that all possible densities of each input $X_k$ have the same support $\mathcal{X}_k$. The main objective is to choose $\bar{f}_k$ as \enquote{close} as possible to the set all potential densities $\mathcal{F}_k$, while accommodating the probability distribution over $\mathcal{F}_k$. Three possibilities for this choice are detailed below: the mixture distribution, the Wasserstein barycenter and the Symmetrical Kullback-Leibler barycenter. For this, we consider here and only here, the following generic notations. We designate by $h$ a random one-dimensional density of distribution $\mathbb{H}$, of which the support is denoted $\mathcal{H}$. We also designate by $\mathcal{S}$ the common support of all realizations of $h$.\\

\textbf{The mixture distribution.} We recall that the mixture distribution \cite{everitt1981finite,titterington1985statistical} of $h$ is defined as 
\begin{equation}
\bar{h}_M = \mathbb{E}_{\mathbb{H}} \left[ h \right] =  \int_{\mathcal{H}} h \; \mathrm{d} \mathbb{H}(h).
\label{fmel}
\end{equation}
In particular, when $\mathbb{H}$ is discrete with support $\mathcal{H} = \lbrace h_{1} , \ldots , h_{m} \rbrace$, we obtain $\bar{h}_M = \sum_{r = 1}^{m} h_{r} \; \mathbb{H}(h_r)$. Moreover, if the support $\mathcal{H}$ is parameterizable, i.e. $\mathcal{H} = \left\lbrace h_{\theta} , \theta \in \Theta \right\rbrace$, we have $\bar{h}_M = \int_{\Theta} h_{\theta}  \pi(\theta) \, \mathrm{d} \theta$, where $\pi$ is the distribution of the parameter $\theta$.\\ 

\textbf{The Symmetrical Kullback-Leibler barycenter.} This barycenter is computed with respect to the so-called \textit{Symmetrical Kullback-Leibler divergence}, which is obtained by symmetrizing the usual Kullback-Leibler divergence \cite{KullbackLeibler1951}. It is defined for two probability measures $\mathbb{Q}_1$ and $\mathbb{Q}_2$ as
\begin{equation}
D_{\SKL} (\mathbb{Q}_1, \mathbb{Q}_2) = \frac{1}{2} \left[ \KL (\mathbb{Q}_1 \vert\vert \mathbb{Q}_2) + \KL (\mathbb{Q}_2 \vert\vert \mathbb{Q}_1) \right], 
\end{equation}
where $\KL( \mathbb{Q}_1 \vert \vert \mathbb{Q}_2 ) = \mathbb{E}_{\mathbb{Q}_1} \left[ \log( \mathrm{d} \mathbb{Q}_1 / \mathrm{d} \mathbb{Q}_2) \right]$ , with $\mathrm{d} \mathbb{Q}_1 / \mathrm{d} \mathbb{Q}_2$ refers to the Radon–Nikodym derivative. The explicit formula of the Symmetrical Kullback-Leibler barycenter is unknown. In the specific case where $\mathbb{H}$ is uniform over a finite set $\{ h_1, \ldots, h_m \}$, a good approximation of  this barycenter is shown in \cite{veldhuis2002centroid} and given by  
\begin{equation}
\bar{h}_{K} \simeq \frac{1}{2} \left[ \bar{h} + \frac{\tilde{h}}{\int_{\mathcal{S}} \tilde{h}} \right],
\label{BarKull}
\end{equation}
where $\bar{h} = 1/m \sum_{r=1}^m h_r$ and $\tilde{h} = \prod_{r=1}^m h_r^{1/m}$ are respectively the arithmetic and geometric means of $\{ h_1, \ldots, h_m \}$. In the general case, we propose to approximate the Symmetrical Kullback-Leibler barycenter as 
\begin{equation}
\bar{h}_{K} \simeq \frac{1}{2} \left[ \bar{h}_M + \frac{ e^{\overline{\log(h)}_M}}{\int_{\mathcal{S}} e^{\overline{\log(h)}_M}} \right],
\label{BarKullcon}
\end{equation}
where $\bar{h}_M$ and $\overline{\log(h)}_M$ are respectively the mixture distributions of the random functions $h$ and $\log(h)$.\\

\textbf{The Wasserstein barycenter distribution.} We remind that the Wassertein distance \cite{givens1984class, villani2003topics} between two distributions $\mathbb{Q}_1$ and $\mathbb{Q}_2$ is defined as 
\begin{equation}
\text{W} (\mathbb{Q}_1, \mathbb{Q}_2) = \inf_{\gamma \in \Gamma(\mathbb{Q}_1,\mathbb{Q}_2)} \left( \mathbb{E}_{\gamma} \left[ ( X - Y )^2 \right] \right)^{1/2},
\end{equation}
where $\Gamma(\mathbb{Q}_1,\mathbb{Q}_2)$ is the set of probabilities of $(X,Y)$ with marginals $\mathbb{Q}_1$ and $\mathbb{Q}_2$. The quantile function of the Wasserstein barycenter \cite{agueh2011barycenters} of a finite uniformly weighted set $\{h_1, \ldots, h_m \}$ is defined as 
\begin{equation}
\bar{q}_{\text{W}} = \frac{1}{m} \sum_{r = 1}^m q_r,
\label{Barwa}
\end{equation}  
where $q_r$ is the quantile function associated to $h_r$. In the general case, we extend Equation \eqref{Barwa} when $h$ is generated according to a distribution $\mathbb{H}$ as 
\begin{equation}
\bar{q}_{\text{W}} = \mathbb{E}_{\mathbb{H}} [q_h],
\label{Barwag}
\end{equation}
where $q_h$ is the quantile function associated to $h$.

\section{Application of GSA2 methodology}
\label{Application}

First, the performance of our methodology is studied through simulated data. More specifically, the drawing density options presented in Section \ref{unique-law} are studied and compared. Moreover, we shed light on the benefit of this approach compared to the \enquote{double-loop} one. Secondly, the methodology is applied on a nuclear case study simulating a severe nuclear reactor accident.

\subsection{Analytical example}
\label{Analytical-example}

To assess the efficiency of the \enquote{single-loop} methodology, we consider the analytical model $\mathcal{M}$ defined in Equation \eqref{Ishimdf}. The inputs $X_1$, $X_2$ and $X_3$ are assumed to be independent. Moreover, their probability distributions $\mathbb{P}_{X_1}$, $\mathbb{P}_{X_2}$ and $\mathbb{P}_{X_3}$ can equiprobably be $\mathbb{P}_U$, $\mathbb{P}_T$ or $\mathbb{P}_N$, where $\mathbb{P}_U $ is the uniform distribution on $[0,1]$, $\mathbb{P}_T$ is the triangular distribution on $[0,1]$ with mode $0.4$, and $\mathbb{P}_N$ is the truncated normal distribution on $[0,1]$ with mean $0.6$ and standard deviation $0.2$.

\smallskip
In practice, this configuration may occurs when for example three experts agree on the input variation ranges but, have different opinions on the nature of the probability distribution. More precisely:
\begin{itemize}
\item the first expert claims that except the range of variation, no other information can be assumed on the uncertain variable;
\item the second adds that the most likely value is 0.4; 
\item the third thinks that the mean and the standard deviation can respectively be assumed equal to 0.6 and 0.2.
\end{itemize}
According to the principle of maximum entropy for expert elicitation \cite{meyer2001eliciting, o2006uncertain}, the information provided by these experts are respectively modeled by the distributions $\mathbb{P}_U$, $\mathbb{P}_T$ and $\mathbb{P}_N$. By assigning equal importance to these three opinions, the uniform distribution on the set $\{\mathbb{P}_U, \mathbb{P}_T, \mathbb{P}_N \}$ seems to be here the most reasonable choice for the second-level uncertainty.   

\smallskip
As a first step, we well approximate the theoretical values of 2$^{\text{nd}}$-level GSA indices for the model $\mathcal{M}$. Subsequently, to study the convergence rates of \enquote{single-loop} estimators, we apply Algorithm \ref{AlgoGSA2} for different sample sizes. To define 1$^{\text{st}}$ and 2$^{\text{nd}}$-level HSIC measures, we use Standardized Gaussian kernel (see Remark \ref{rem22}) for all vector quantities and the kernels presented in Section \ref{RKHS-choices} for non-vectorial quantities of interest.

\subsubsection{Computation of theorical values}
\label{Theorical-values}
We focus here on Option \ref{QI-R2} of Algorithm \ref{AlgoGSA2}, the other quantities of interest are studied in Section \ref{GSA2-Other-QI}. To approximate the theoretical values of 2$^{\text{nd}}$-level HSIC measures and indices, we consider the set of the $n_1 = 27$ possible 3-tuples of input probability distributions. For each input distribution, the 1$^{\text{st}}$-level HSIC measures and indices are computed using a sample of size $n_2 = 1000$, generated according to the prior input density. The theoretical 2$^{\text{nd}}$-level HSIC measures are computed: 
$$\HSIC (\mathbb{P}_{X_1} , \mathcal{R} ) = 0.0414, \; \HSIC (\mathbb{P}_{X_2} , \mathcal{R} ) = 0.0261 \mbox{ and } \HSIC (\mathbb{P}_{X_3} , \mathcal{R} ) = 0.0009.$$
The theoretical 2$^{\text{nd}}$-level HSIC indices are also computed:
$$ \R2^2_{\HSIC} (\mathbb{P}_{X_1} , \mathcal{R} ) = 0.4152, \; \R2^2_{\HSIC} (\mathbb{P}_{X_2} , \mathcal{R} ) = 0.2516 \mbox{ and }  \R2^2_{\HSIC} (\mathbb{P}_{X_3} , \mathcal{R} ) = 0.0086.$$

In this example, we observe that $\R2^2_{\HSIC} (\mathbb{P}_{X_1} , \mathcal{R})$ is significantly larger than the other two indices, while $\R2^2_{\HSIC} (\mathbb{P}_{X_3} , \mathcal{R})$ is negligible. Based on these results, the lack of knowledge on $\mathbb{P}_{X_3}$ is not responsible for the variability of 1$^{\text{st}}$-level HSIC indices. This distribution can simply be set to a reference one. Furthermore, the impact of $\mathbb{P}_{X_1}$ uncertainty is by far the largest and the one of $\mathbb{P}_{X_2}$ remains non-negligible. Therefore, characterization efforts should be targeted in priority on $\mathbb{P}_{X_1}$, followed-up by $\mathbb{P}_{X_2}$.

\subsubsection{GSA2 with our single-loop approach}
\label{Test&options}

In the sequel, 2$^{\text{nd}}$-level HSIC estimators using the mixture distribution, the Wasserstein barycenter and the Symmetrical Kullback-Leibler barycenters are respectively denoted $\widetilde{\HSIC}_M (\mathbb{P}_{X_k}, \mathcal{R})$, $\widetilde{\HSIC}_W (\mathbb{P}_{X_k}, \mathcal{R})$ and $\widetilde{\HSIC}_K (\mathbb{P}_{X_k}, \mathcal{R})$. Similarly, the 2$^{\text{nd}}$-level indices are denoted  
$\widetilde{\R2}^2_{\HSIC,\M} (\mathbb{P}_{X_k} , \mathcal{R})$, $\widetilde{\R2}^2_{\HSIC,\WW} (\mathbb{P}_{X_k} , \mathcal{R})$ and $\widetilde{\R2}^2_{\HSIC,\K} (\mathbb{P}_{X_k} , \mathcal{R})$.

\begin{figure}[h!]
\centering
\includegraphics[scale=0.15]{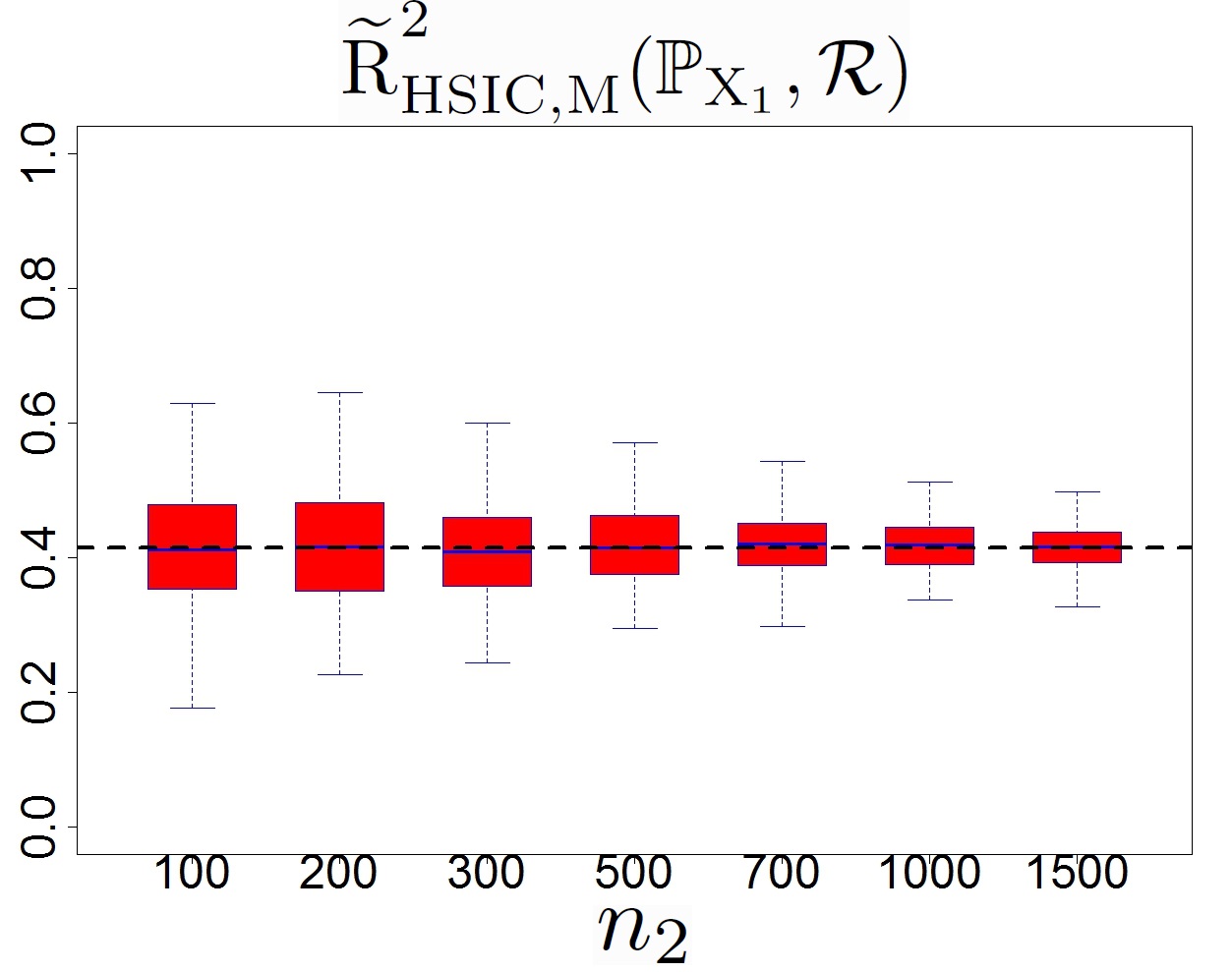} 
\includegraphics[scale=0.15]{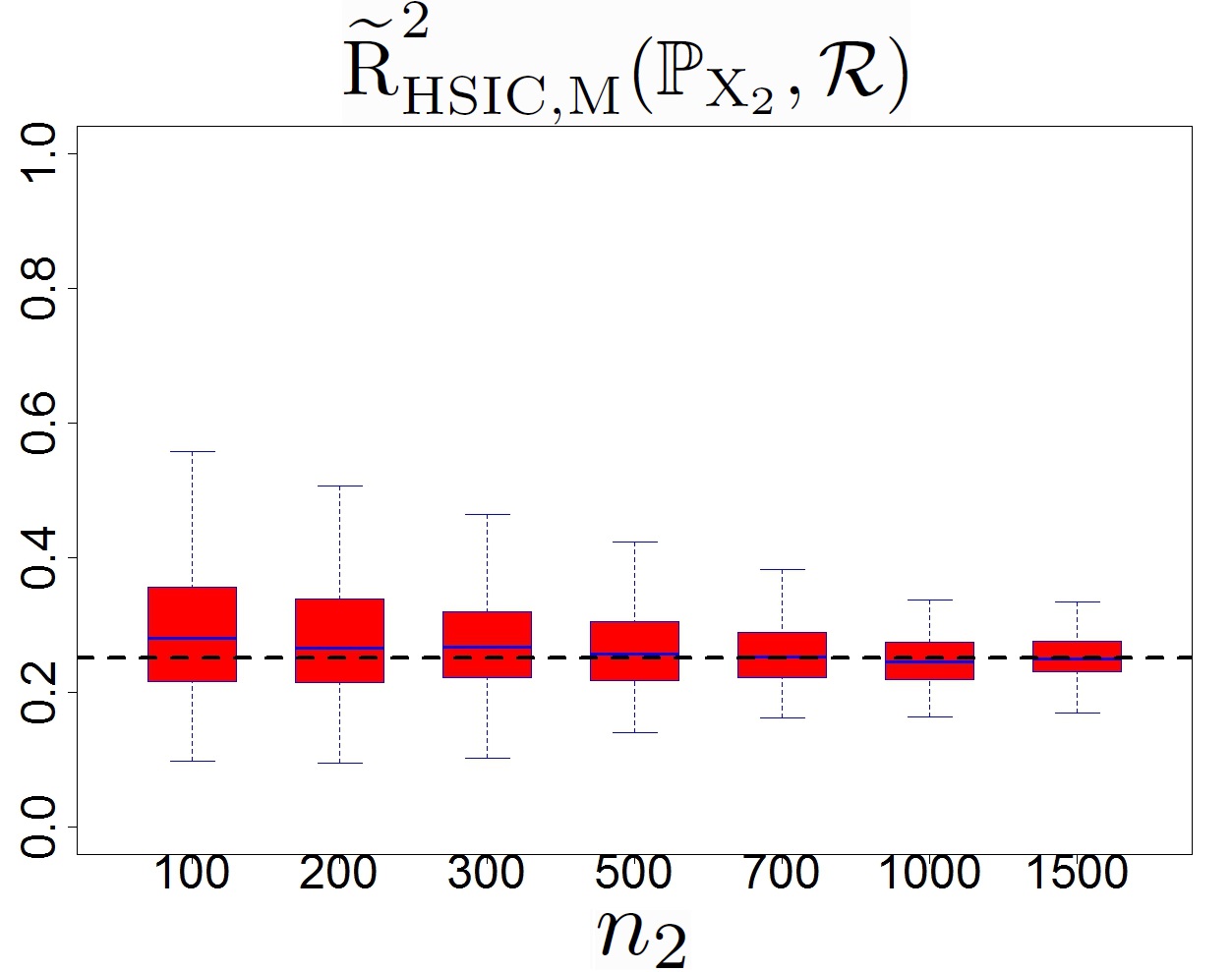}
\includegraphics[scale=0.15]{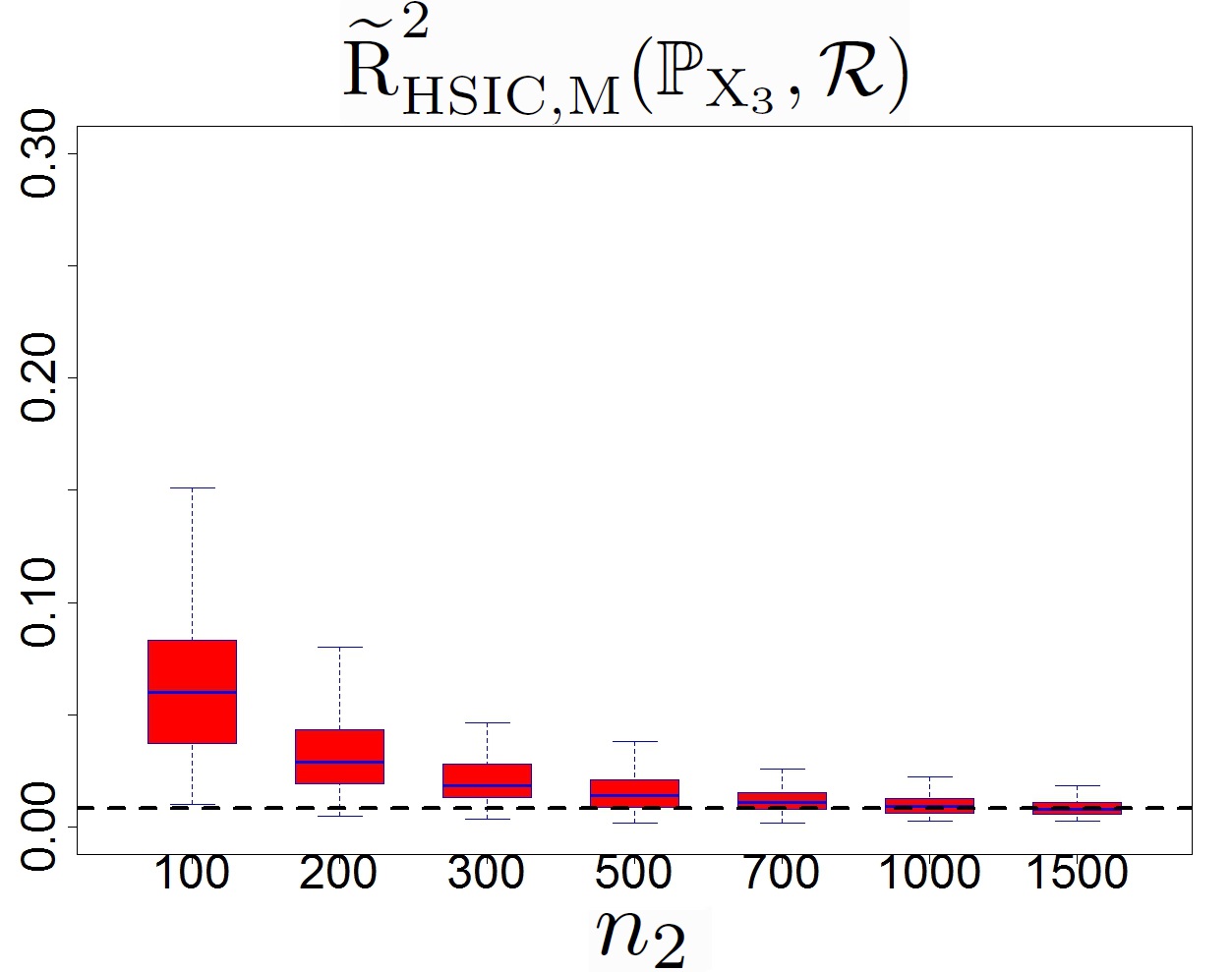}
\centering
\vspace{0.05cm}
\includegraphics[scale=0.15]{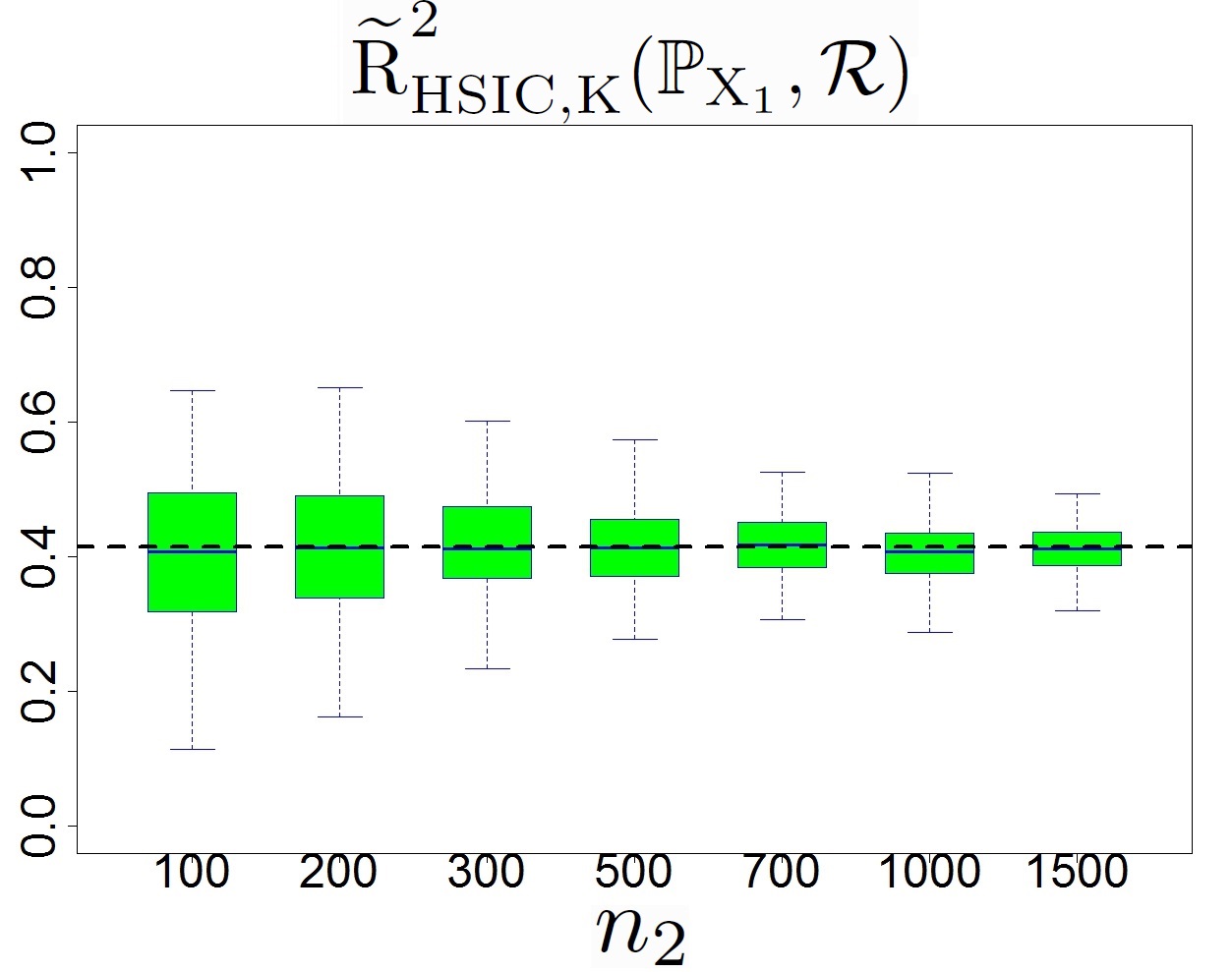} 
\includegraphics[scale=0.15]{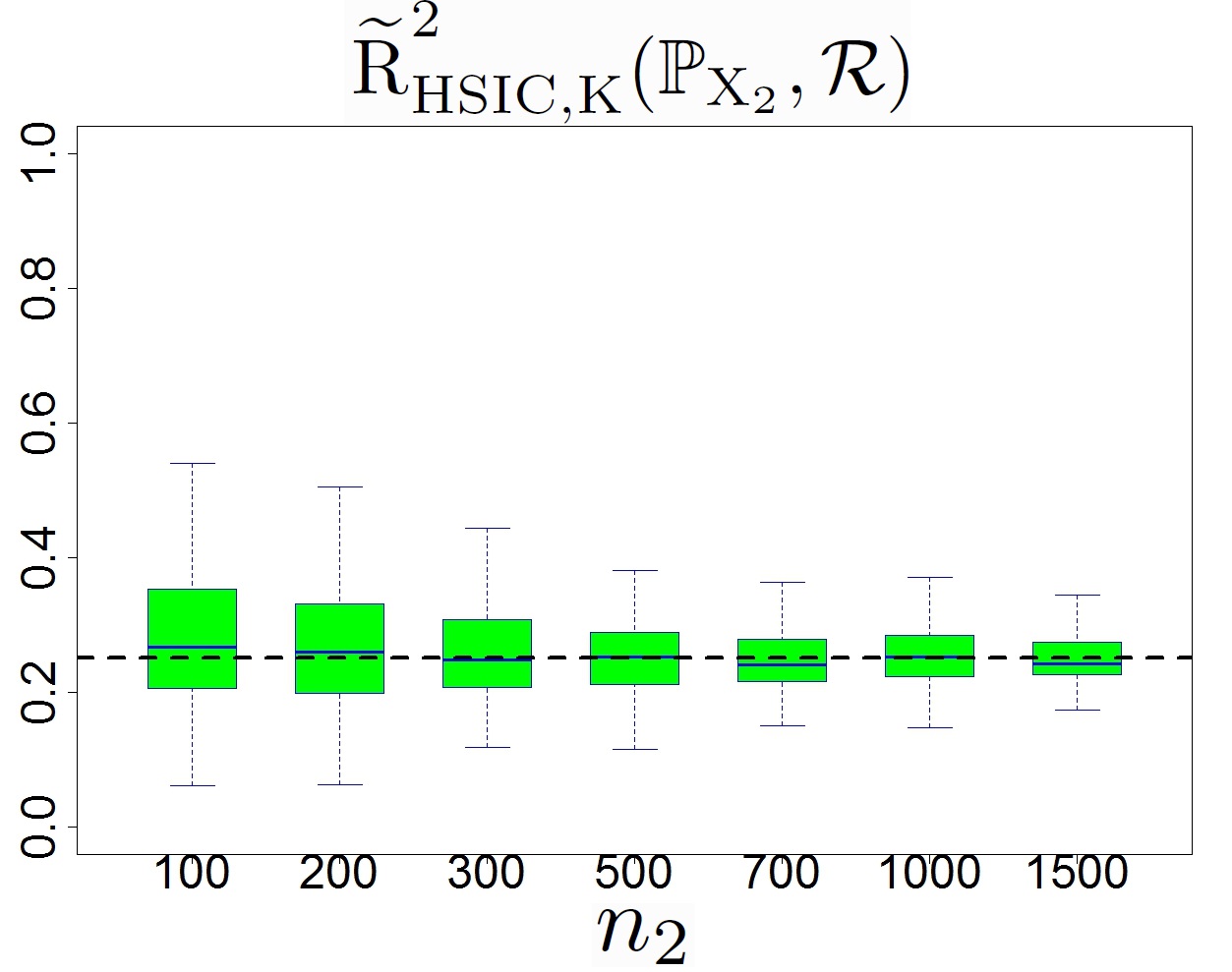}
\includegraphics[scale=0.15]{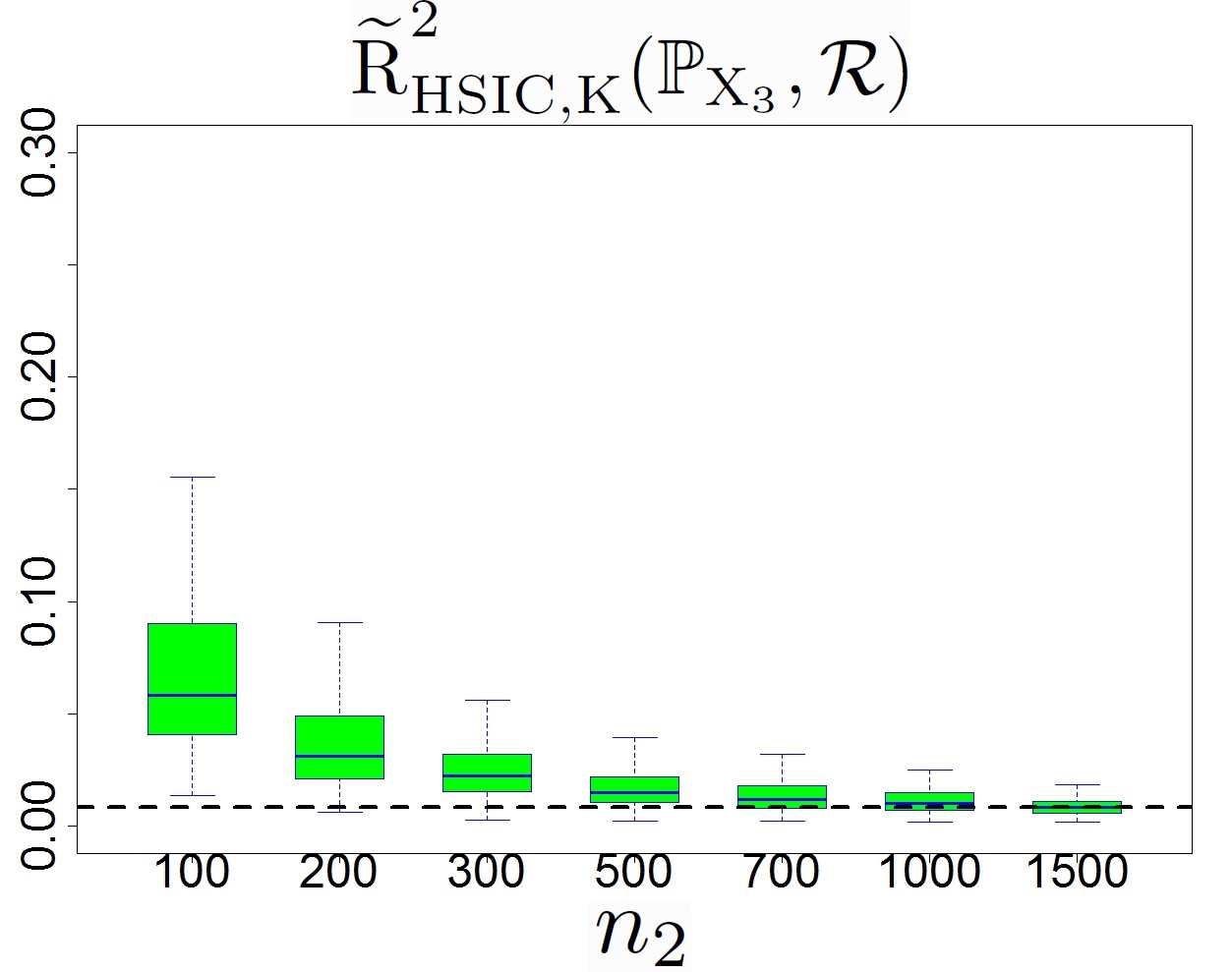}
\centering
\vspace{0.05cm}
\includegraphics[scale=0.15]{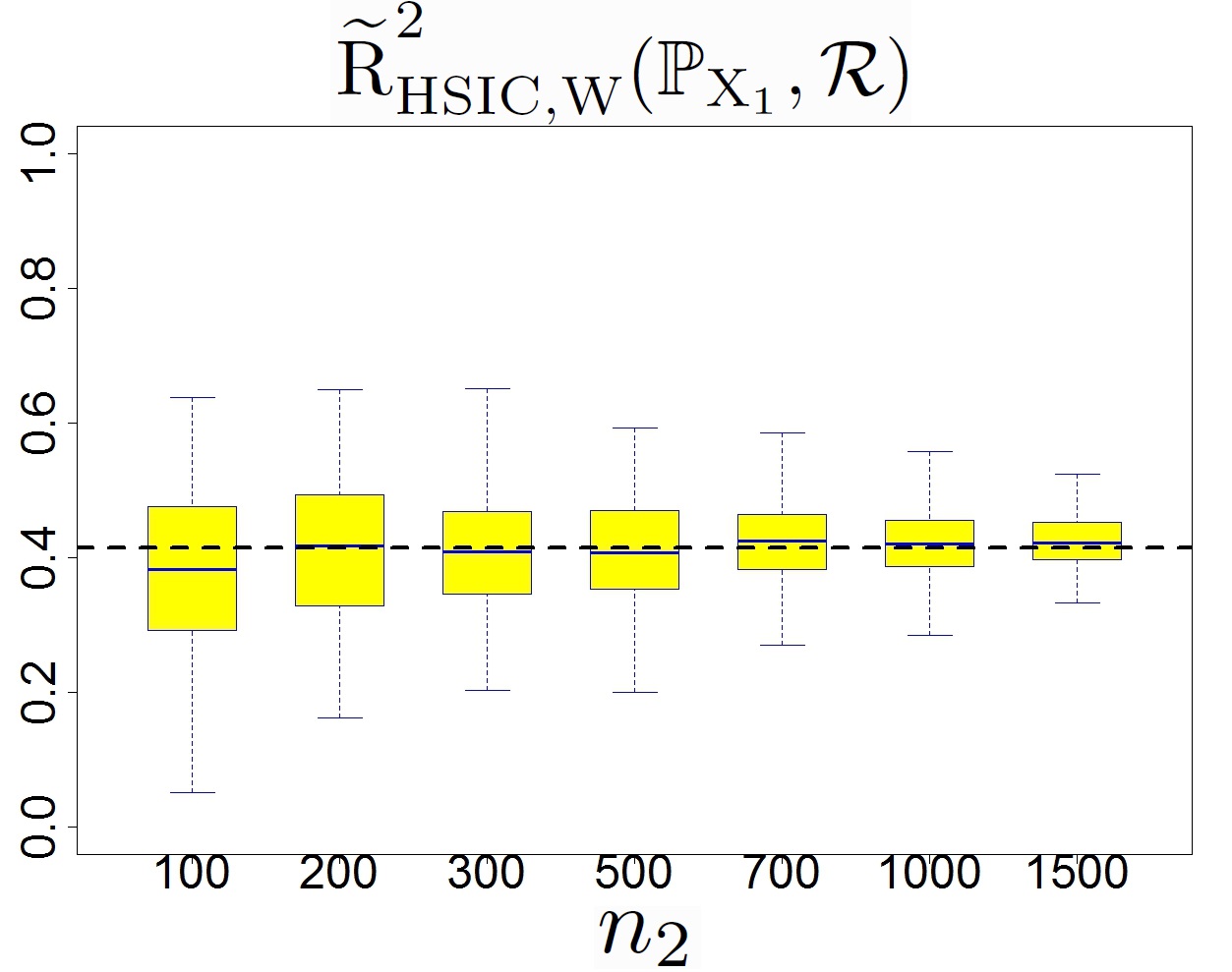} 
\includegraphics[scale=0.15]{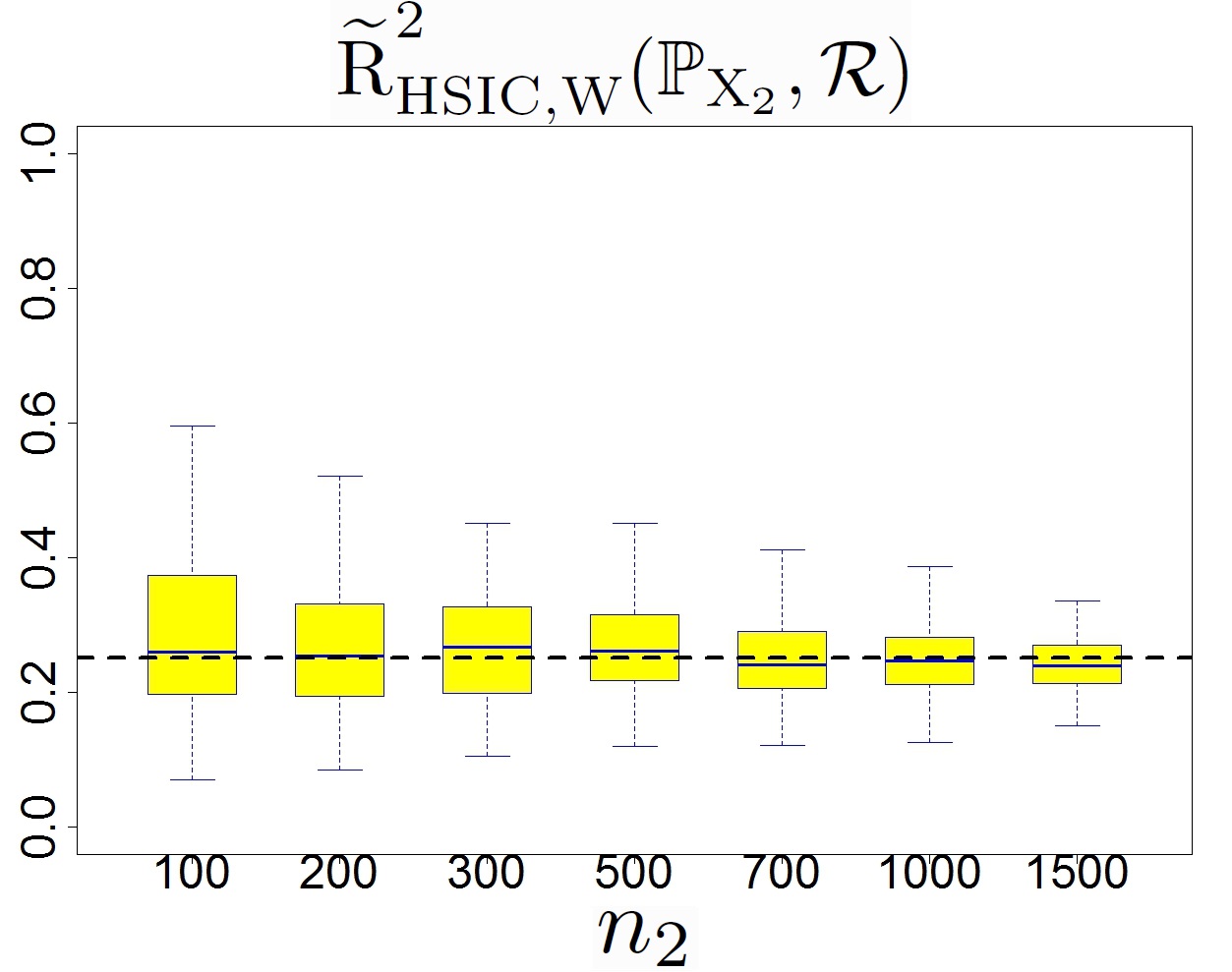}
\includegraphics[scale=0.15]{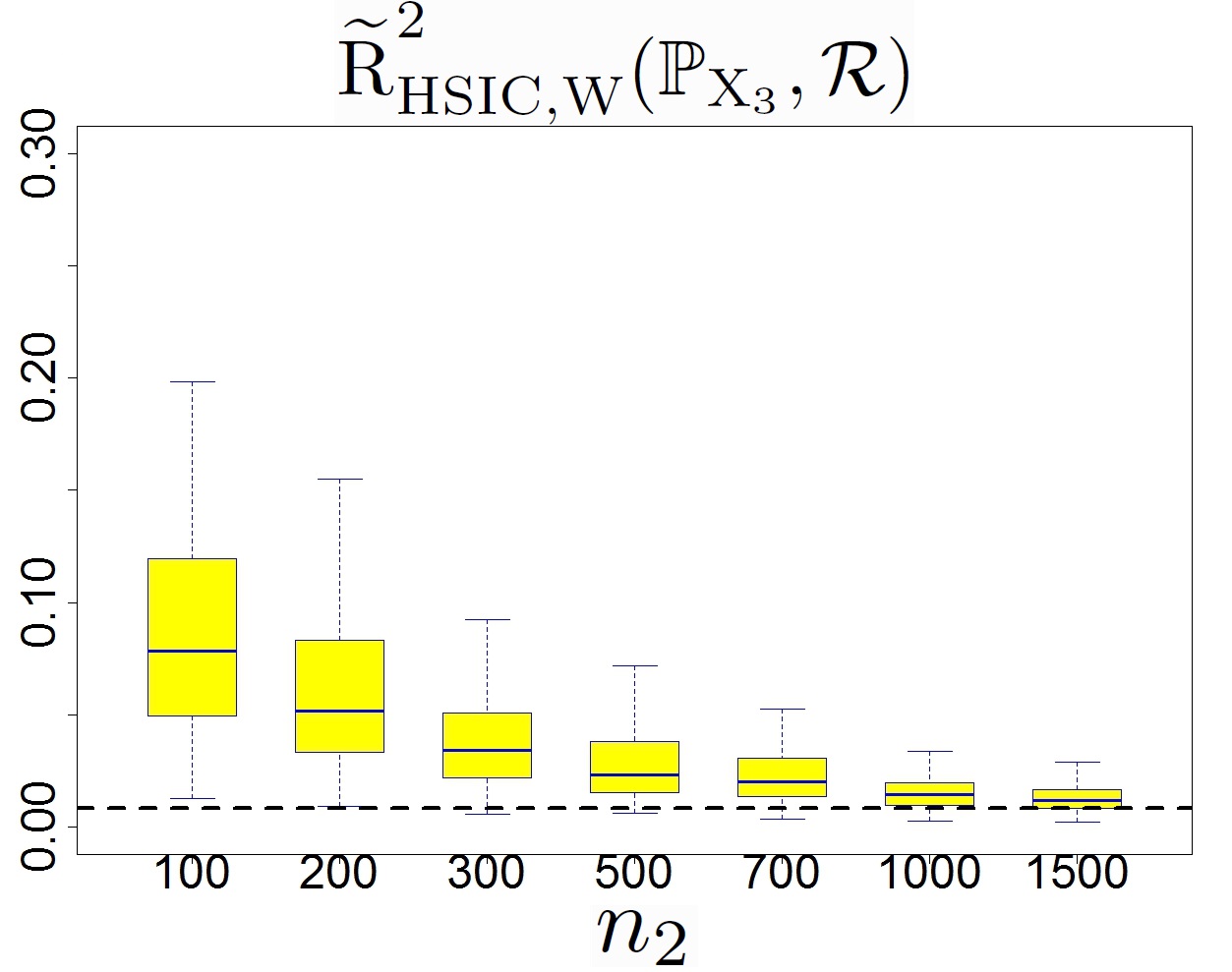}
\caption{Model $\mathcal{M}$ -- Convergence plots of the estimators $\widetilde{\R2}^2_{\HSIC,\M} (\mathbb{P}_{X_k} , \mathcal{R})$, $\widetilde{\R2}^2_{\HSIC,\WW} (\mathbb{P}_{X_k} , \mathcal{R})$ and $\widetilde{\R2}^2_{\HSIC,\K} (\mathbb{P}_{X_k} , \mathcal{R})$, with respect to the sample size $n_2$. Theoretical values are represented in dotted lines.}
\label{HSIC23es} 
\end{figure}

\smallskip
To study the convergence rate of the \enquote{single-loop} estimators, we apply Algorithm \ref{AlgoGSA2} from samples with sizes ranging from $n_2 = 100 $ to $ n_2 = 1500$. For each sample size, the estimations are repeated independently 200 times using independent samples. Results are given by Figure \ref{HSIC23es}, where the theoretical values of $\R2^2_{\HSIC} (\mathbb{P}_{X_k} , \mathcal{R})$ are represented in dotted lines. Visually,  the estimators based on the mixture distribution and the Symmetrical Kullback-Leibler barycenter seem to perform similarly both for small and large sample sizes. In particular, the dispersion of these estimators are satisfying from $n_2=700$. In contrast, the Wasserstein barycenter estimators are less accurate (higher dispersion) compared to the previous estimators, especially for small and medium size samples (i.e. $n_2$ in $[300,700]$).

\smallskip
A more pragmatic way to compare these drawing densities, is to compare the estimators ability to correctly rank the input distributions based on their influence level. To do so, we compute for each sample size, the percentage of times the theoretical ranking and the one given by the estimators match. The results are presented in Table \ref{tauxbon2}. As expected, the estimators based on the mixture distribution and the Symmetrical Kullback-Leibler barycenter outperform those based on the Wasserstein barycenter, and this regardless of the sample size. This can be explained by the fact that the likelihood ratio $f/f_W$ is very high in the neighborhoods of 0 and 1. Furthermore, the Kullback-Leibler barycenter seems to give slightly better results for small samples $n_2 \leq 300$, the reverse is true from $n_2 = 500$.

\begin{table}[!h]
\begin{center}
\begin{tabular}{|D{2.9cm}|D{1.3cm}|D{1.3cm}|D{1.3cm}|D{1.3cm}|D{1.3cm}|D{1.3cm}|D{1.3cm}|}\hline
$n_2$ & $100$ & $200$ & $300$ &$500$ & $700$ &$1000$ & $1500$ \tabularnewline\hline \hline
$\widetilde{\R2}^2_{\HSIC,\M} (\mathbb{P}_{X_k} , \mathcal{R})$ &$74\%$ & $79\%$ & $84\%$& $94.5\%$ & $97\%$ & $100\%$ & $100 \%$\tabularnewline\hline
$\widetilde{\R2}^2_{\HSIC,\K} (\mathbb{P}_{X_k} , \mathcal{R})$ &$75.5\%$ & $79\%$ & $87\%$& $92\%$ & $97\%$ & $99.5\%$ & $99.5 \%$\tabularnewline\hline
$\widetilde{\R2}^2_{\HSIC,\WW} (\mathbb{P}_{X_k} , \mathcal{R})$ &$57.5 \%$ & $71\%$ & $77\%$& $82\%$ & $91\%$ & $93.5\%$ & $98\%$\tabularnewline\hline
\end{tabular} 
\caption{Model $\mathcal{M}$ -- Good ranking rates of $(\mathbb{P}_{X_1}, \mathbb{P}_{X_2}, \mathbb{P}_{X_3}$) using the estimators $\widetilde{\R2}^2_{\HSIC,\M} (\mathbb{P}_{X_k} , \mathcal{R})$, $\widetilde{\R2}^2_{\HSIC,\K} (\mathbb{P}_{X_k} , \mathcal{R})$ and $\widetilde{\R2}^2_{\HSIC,\WW} (\mathbb{P}_{X_k} , \mathcal{R})$, with respect to the sample size $n_2$.}
\label{tauxbon2}
\end{center}
\end{table}

\subsubsection{Comparison with Monte Carlo \enquote{double-loop} approach}
\label{DoubleVsSingle}

We compare now the performance of the \enquote{single-loop} and \enquote{double-loop} approaches, in terms of convergence rates of estimators. To do so, we consider a total simulation budget of $n = 1026$ for both approaches. More precisely for the \enquote{double-loop} approach, a sample of size $n_2 = 38$ is generated for each possible 3-tuple of input distributions (for a total number of $n = n_1 \times n_2 = 1026$ simulations). The associated estimators are denoted $\widehat{\R2}^2_{\HSIC} (\mathbb{P}_{X_k}, \mathcal{R})$ with $k \in \{ 1, 2, 3 \}$. Concerning the \enquote{single-loop} approach, we apply Algorithm \ref{AlgoGSA2} with $n_2 = 1026$ and we compute the estimators $\widetilde{\R2}^2_{\HSIC,\M} (\mathbb{P}_{X_k} , \mathcal{R})$ and $\widetilde{\R2}^2_{\HSIC,\K} (\mathbb{P}_{X_k} , \mathcal{R})$ with $k \in \{ 1, 2, 3 \}$.

\begin{figure}[!h]
\centering
\includegraphics[scale=0.26]{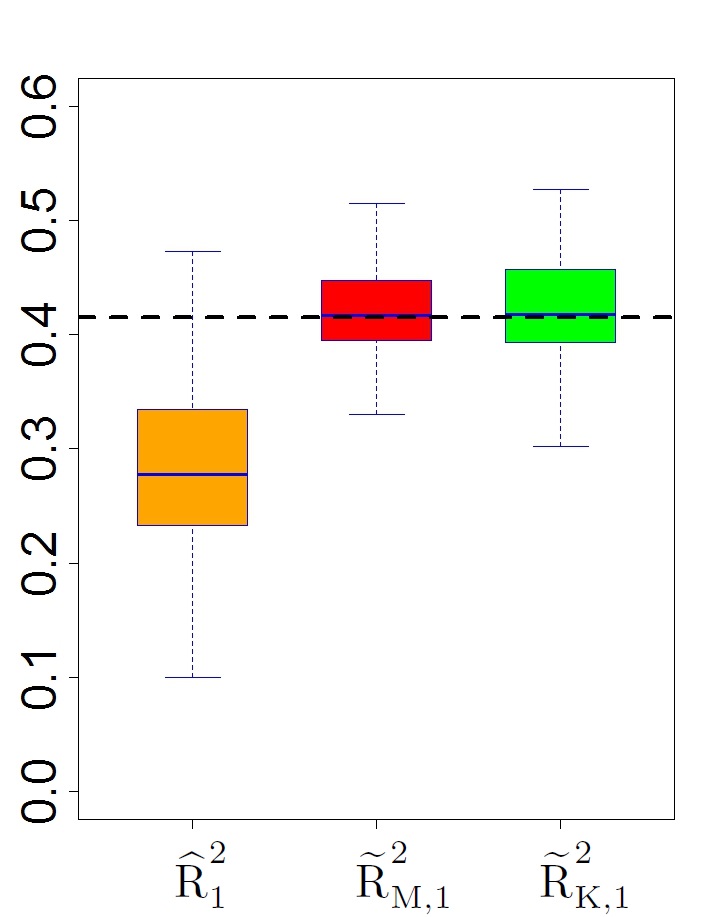} 
\includegraphics[scale=0.26]{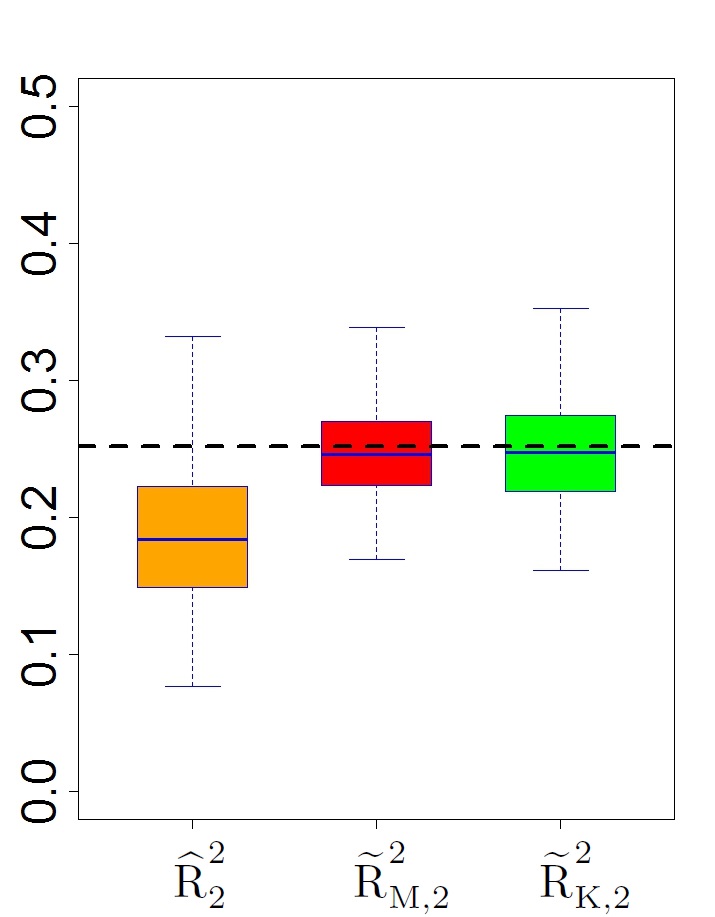}
\includegraphics[scale=0.26]{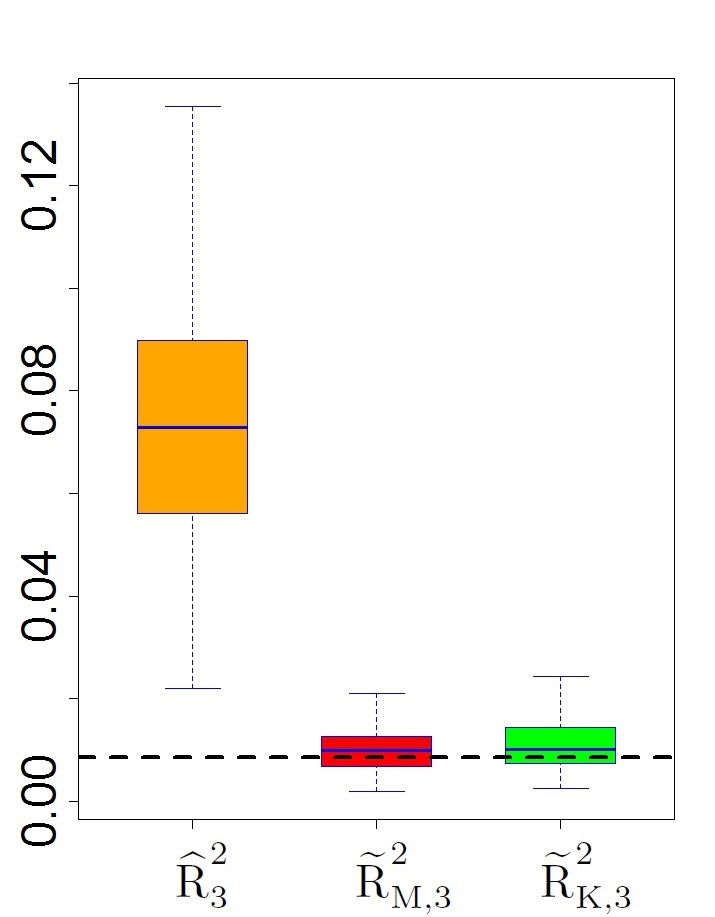}
\caption{Model $\mathcal{M}$ -- Comparison of 2$^{\text{nd}}$-level HSIC indices estimated by the \enquote{single-loop} and \enquote{double-loop} approaches for $n = 1026$. The estimators are denoted, $\widehat{\R2}^2_{k}$ for $ \widehat{\R2}^2_{\HSIC} (\mathbb{P}_{X_k}, \mathcal{R})$, $\widetilde{\R2}^2_{\M,k} $ for $ \widetilde{\R2}^2_{\HSIC,\M} (\mathbb{P}_{X_k} , \mathcal{R})$ and $\widetilde{\R2}^2_{\K,k} $ for $ \widetilde{\R2}^2_{\HSIC,\K} (\mathbb{P}_{X_k} , \mathcal{R})$. Theoretical values are represented in dotted lines.}
\label{HSIC2MWK}
\end{figure}

\smallskip
Each estimation is repeated 200 times with independent Monte Carlo samples. The estimator boxplots are shown by Figure~\ref{HSIC2MWK}, where the theoretical values are represented in dotted lines. The \enquote{double-loop}  estimators show much more variability than the \enquote{single-loop} ones. Also, notice that the \enquote{single-loop} estimators are much less biased than the \enquote{double-loop} ones. Our approach significantly outperforms the \enquote{double loop}. This conclusion is also supported by the \textit{good ranking rates} presented in Table \ref{tauxbonDB}. 
A reasonable explanation for the benefit of the \enquote{single-loop} approach, may be the simulation budget for 1$^{\text{st}}$-level HSIC. Indeed, given a total budget of $n$ simulations, each 1$^{\text{st}}$-level HSIC is computed using $n_2 = n$ for the \enquote{single-loop} approach, against $n_2 = n/n_1$ for the \enquote{double loop} approach. Although the prior estimators converge faster than the alternative ones, the total simulation number is drastically reduced when using the \enquote{double-loop} approach. 

\smallskip
For this same model $\mathcal{M}$, other numerical studies with different
hypothesis on input distribution uncertainty have been performed and yield similar results and conclusions.

\begin{table}[!h]
\begin{center}
\begin{tabular}{|D{3.5cm}|D{3.5cm}|D{3.5cm}|}\hline
Double loop  & \multicolumn{2}{c|}{Single loop} \tabularnewline\hline \hline
$\widehat{\R2}^2_{\HSIC} (\mathbb{P}_{X_k}, \mathcal{R})$ & $\widetilde{\R2}^2_{\HSIC,\M} (\mathbb{P}_{X_k} , \mathcal{R})$ & $\widetilde{\R2}^2_{\HSIC,\K} (\mathbb{P}_{X_k} , \mathcal{R})$ \tabularnewline\hline 
$67.5 \%$ & $100 \%$ & $99 \%$  \tabularnewline\hline
\end{tabular} 
\caption{Model $\mathcal{M}$ -- Comparison of good ranking rates of \enquote{single-loop} and \enquote{double-loop} estimators, for $n = 1026$.}
\label{tauxbonDB}
\end{center}
\end{table}

\subsubsection{GSA2 with other quantities of interest}
\label{GSA2-Other-QI}

It is fair to wonder whether  GSA2 conclusions vary if we decide to choose other quantities of interest. To answer that, we deal with Options \ref{QI-RK}, \ref{Pval-asy} and \ref{Pval-perm} of Algorithm \ref{AlgoGSA2}. In all cases, we keep the same kernel choices as described at the beginning of Section \ref{Analytical-example}. Let us examine these possibilities one-by-one.\\

$\bullet$ \textbf{Ranking by $\R2^2_{\HSIC}$.} Before looking closely at the simulation results, one can notice that the convergence of 1$^{\text{st}}$-level $\R2^2_{\HSIC}$ estimators systematically implies the convergence of those by ranking. Therefore, the estimators of GSA2 indices of Option \ref{QI-RK} converge faster than those of Option \ref{QI-R2}. Moreover, according to Section \ref{Test&options} results, a sample of size $n_2 = 1000$ of the drawing density is sufficient to accurately estimate the indices. We thus obtain:
$$ \R2^2_{\HSIC} (\mathbb{P}_{X_1} , \mathcal{R} ) = 0.3830, \; \R2^2_{\HSIC} (\mathbb{P}_{X_2} , \mathcal{R} ) = 0.0958 \mbox{ and }  \R2^2_{\HSIC} (\mathbb{P}_{X_3} , \mathcal{R} ) \simeq 0. $$

The gaps of theses values are more meaningful compared to those presented in Section \ref{Theorical-values}. This is likely related to the stability of the ranking compared to GSA1 indices. Indeed, only significant variations of GSA1 indices contribute to GSA2 indices using the ranking. We safely conclude that $\mathbb{P}_{X_1}$ is the main contributor for the ranking uncertainty; less characterization efforts are required.

\medskip
$\bullet$ \textbf{P-values vector.} When considering Options \ref{Pval-asy} or \ref{Pval-perm} as the quantity of interest, two points are highlighted. Firstly, the estimators of GSA2 indices show a large variance, regardless of the p-value estimation method (Gamma approximation or permutations), even for very large $n_2$ such as $n_2 = 5000$. In addition, the three estimated GSA2 indices are small (not exceeding 0.2). To help understanding these results, we focus on the estimated p-values for each possible input distribution. To do so, we use the permutation method with $B=1000$ resamplings. The results show that the p-values associated to $X_1$ and $X_2$ are almost equal to zero (exactly zero numerically), regardless of the input distribution. Moreover, the p-values associated to $X_3$ are very low and in most cases below $10^{-5}$. Therefore, the high variance of GSA2 indices is due to the difficulty of accurately estimating each p-value. In this case, Options \ref{Pval-asy} and \ref{Pval-perm} are not relevant: the independence hypothesis is not reliable and this, regardless of the input distribution. The 2$^{\text{nd}}$-level input uncertainties have no impact on these GSA1 results. 

\subsection{Nuclear safety application}
\label{Nuclear-application}

\begin{figure}[!h]
\centering
\includegraphics[scale=0.26]{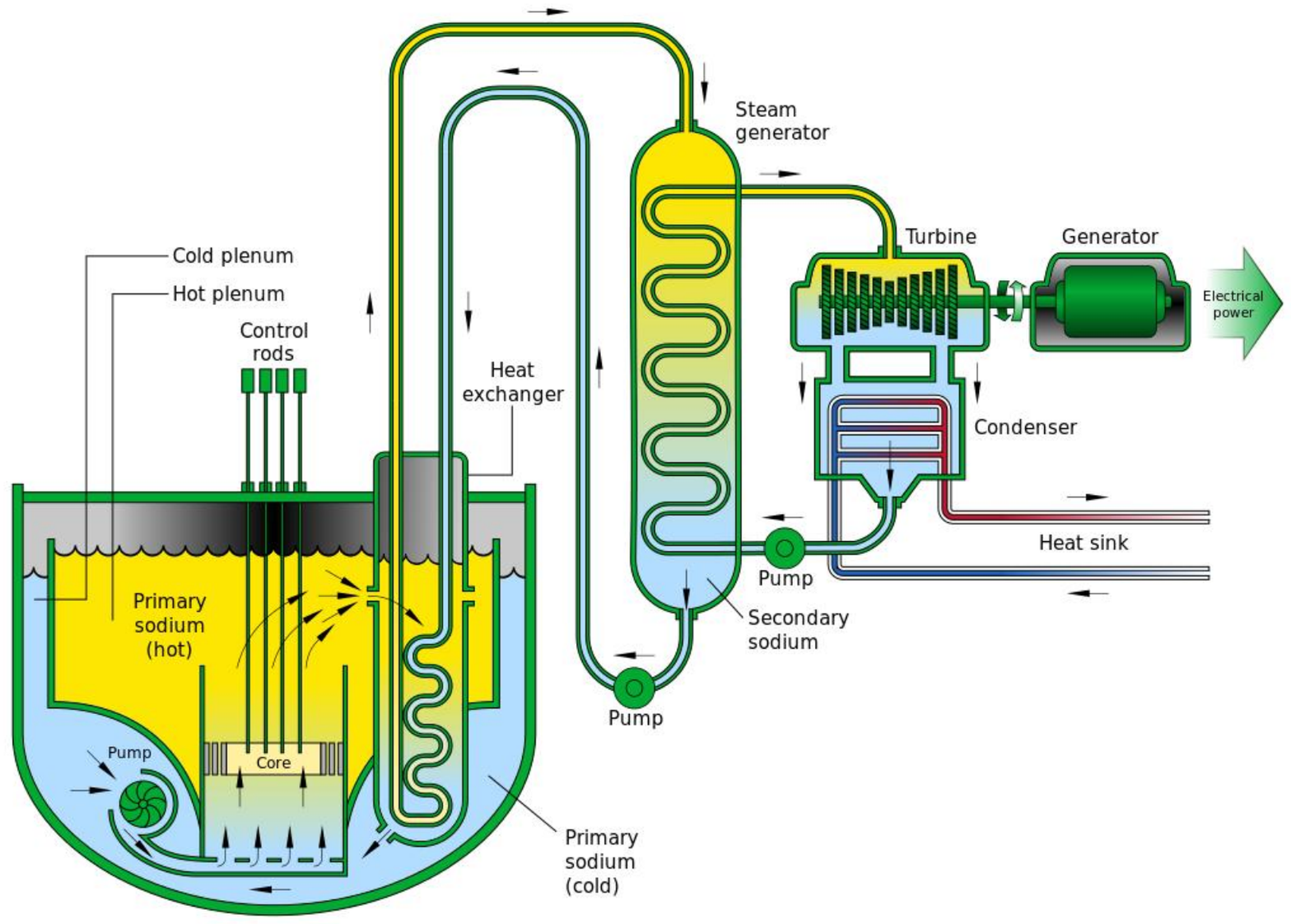} 
\caption[]{MACARENa application -- Basic architecture of a Sodium-cooled Fast Reactor.}
\label{RNR}
\end{figure}

Within the framework of 4$^{\text{th}}$-generation sodium-cooled fast reactor ASTRID: Advanced Sodium Technological Reactor for Industrial Demonstration (see Figure \ref{RNR}), the CEA (French \emph{Commissariat à l'Énergie atomique et aux Énergies alternatives}) provides numerical tools to model severe accident scenarios  and assess the safety. Among them, a numerical tool called MACARENa (French: \emph{Modélisation de l'ACcident d'Arrêt des pompes d'un REacteur refroidi au sodium}) developed by \cite{droin2017physical} simulates a primary phase of an Unprotected Loss Of Flow (ULOF) accident. During this type of accident, the power loss of primary pumps and the dysfunction of shutdown systems cause a gradual decrease of the sodium flow in the primary circuit, which subsequently may increase the temperature of sodium until it boils. This can lead to a degradation of several components and structures of the reactor core. 

\smallskip
Previous GSA studies were performed on MACARENa simulator with several tens of uncertain parameters whose pdf were assumed to be known and set at a reference pdf. These studies show that only 3 input parameters mainly impact the accident transient predicted by MACARENa, namely:
 
\begin{itemize}
\item[$\bullet$] $X_1$: error of measurement on external pressure loss,
\item[$\bullet$] $X_2$: primary half-flow time,
\item[$\bullet$] $X_3$: Lockart-Martinelli correction value.
\end{itemize} 

\smallskip 
However, due to lack of data and knowledge, uncertainty remains on the distributions $\mathbb{P}_{X_1}$, $\mathbb{P}_{X_2}$ and $\mathbb{P}_{X_3}$. To take into account this uncertainty, the nature of each input distribution is assumed to be known, but with one uncertain parameter, as described in Table \eqref{LawOfLaw}. The notations $\mathcal{N}_t(a,b,m,\sigma)$, $\T(a,b,c)$ and $\mathcal{U}(a,b)$ are respectively, the truncated normal distribution of mean $m$ and standard deviation $\sigma$ on $[a,b]$, the triangular law on $[a,b]$ with mode $c$ and the uniform  distribution on $[a,b]$. The identification of these uncertainties is based on expert advice. More specifically, the uncertainty on $\sigma$ stems from a prior knowledge (no available data), while the uncertainties on $c$ and $m$ are due to their estimation using few existing partial data.
\begin{table}[!h]
\begin{center}
\begin{tabular}{|D{2.5cm}|D{3.5cm}|D{3.5cm}|}\hline
Law of input & Nature & Uncertain parameter \tabularnewline\hline \hline
$\mathbb{P}_{X_1}$ & $\mathcal{N}_t(-0.1,0.1,0,\sigma)$ & $\sigma \sim \mathcal{U}(0.03,0.05)$
\tabularnewline\hline
$\mathbb{P}_{X_2}$ & $\T(0,20,c)$ & $c \sim \mathcal{U}(8,15)$
\tabularnewline\hline
$\mathbb{P}_{X_3}$ & $\T(0.8,2,m)$  & $m \sim \mathcal{U}(1,1.5)$
\tabularnewline\hline
\end{tabular} 

\smallskip
\caption[]{MACARENa application -- Uncertainties on $\mathbb{P}_{X_1}$,$\mathbb{P}_{X_2}$ and $\mathbb{P}_{X_3}$ distributions.}
\label{LawOfLaw}
\end{center}
\end{table}

\smallskip 
Among the outputs computed by MACARENa simulator to describe the ULOF accident, we focus on the first instant of sodium boiling denoted $Y$. To assess the impact of input distributions on GSA1 results of $Y$, we apply Algorithm \ref{AlgoGSA2} with Option \ref{QI-R2}. We use the mixture density for the unique drawing and the same kernel choices as in Section \ref{Analytical-example}. Moreover, we consider a Monte Carlo sample of size $n_2 = 1000$ for the unique drawing. This choice is motivated by two main reasons. Firstly, the MACARENa simulation cost (between 2 and 3 hours on average) which limits the total number of simulations. Secondly, the analytical three-dimensional example of Section \ref{Analytical-example} for which a budget of 1000 simulations give good results. In addition, we consider a Monte Carlo sample of $n_1 =200$ 3-tuples of input distribution. These two choices for $n_1$ and $n_2$ will numerically be justified later, by studying the stability of estimators.  
Algorithm \ref{AlgoGSA2} gives the following GSA2 indices:
$$\widetilde{\R2}^2_{\HSIC,\M} (\mathbb{P}_{X_1} , \mathcal{R}) = 0.5341, \; \widetilde{\R2}^2_{\HSIC,\M} (\mathbb{P}_{X_2} , \mathcal{R}) = 0.3317 \mbox{ and } \widetilde{\R2}^2_{\HSIC,\M} (\mathbb{P}_{X_3} , \mathcal{R})= 0.0753.$$

\smallskip
Consequently, the uncertainty on $\mathbb{P}_{X_1}$ mainly impacts GSA1 results, followed by $\mathbb{P}_{X_2}$, while the impact of $\mathbb{P}_{X_3}$ is negligible. To improve the robustness of GSA1 results, characterization efforts should then focus primarily on $\mathbb{P}_{X_1}$. 
A deeper analysis of the 200 results of GSA1 shows that the input $X_2$ is always the most predominant. Surprisingly, $X_2$ whose distribution is not the most influential on GSA1 results is the most influential on $Y$.
This example illustrates, if necessary, that the information captured by GSA2 is different but complementary to that of GSA1.

\smallskip
To assess the accuracy of the estimation of GSA2 indices, we use a non-asymptotic bootstrapping approach \cite{efron1994introduction}. For this, we first generate Monte Carlo subsamples with replacement from the initial sample (of 1000 simulations), then we re-estimate 2$^{\text{nd}}$-level $\R2^2_{\HSIC}$ using these samples. More specifically, we consider subsamples of sizes $n_2 = 100$ to $n_2 = 800$. For each size, the estimation is repeated independently $B=20$ times. Furthermore, to reduce computational efforts, we consider a sample of distributions of reduced size $n_1 = 30$ and generated with a space-filling approach. More precisely, the vector $(\sigma,c,m)$ is sampled with a Maximum Projection Latin Hypercube Design \cite{joseph2015maximum} of size $n_1 = 30$ and defined on the cubic domain $[0.03,0.05] \times [8,15] \times [1,1.5]$.

\begin{figure}[ht]
\centering
\includegraphics[scale=0.16]{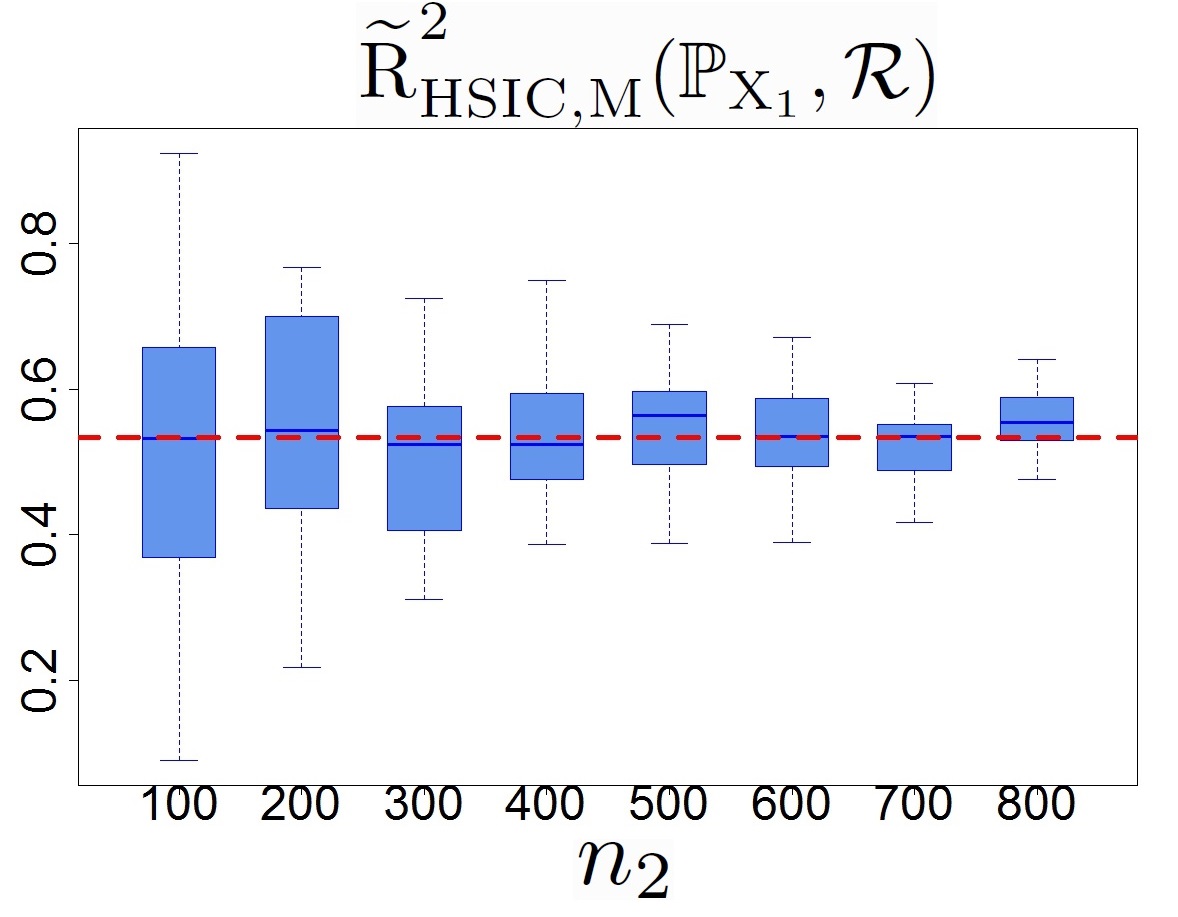} 
\includegraphics[scale=0.16]{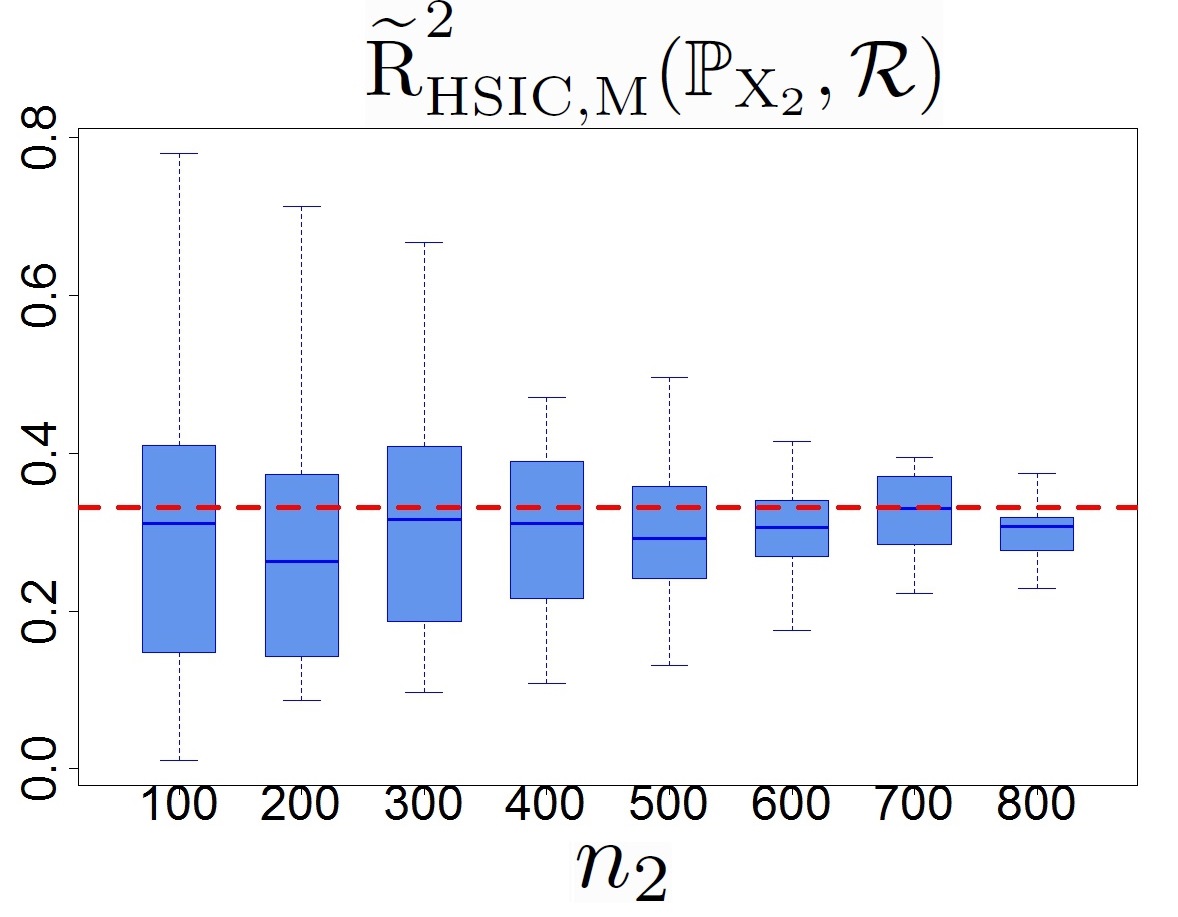}
\includegraphics[scale=0.16]{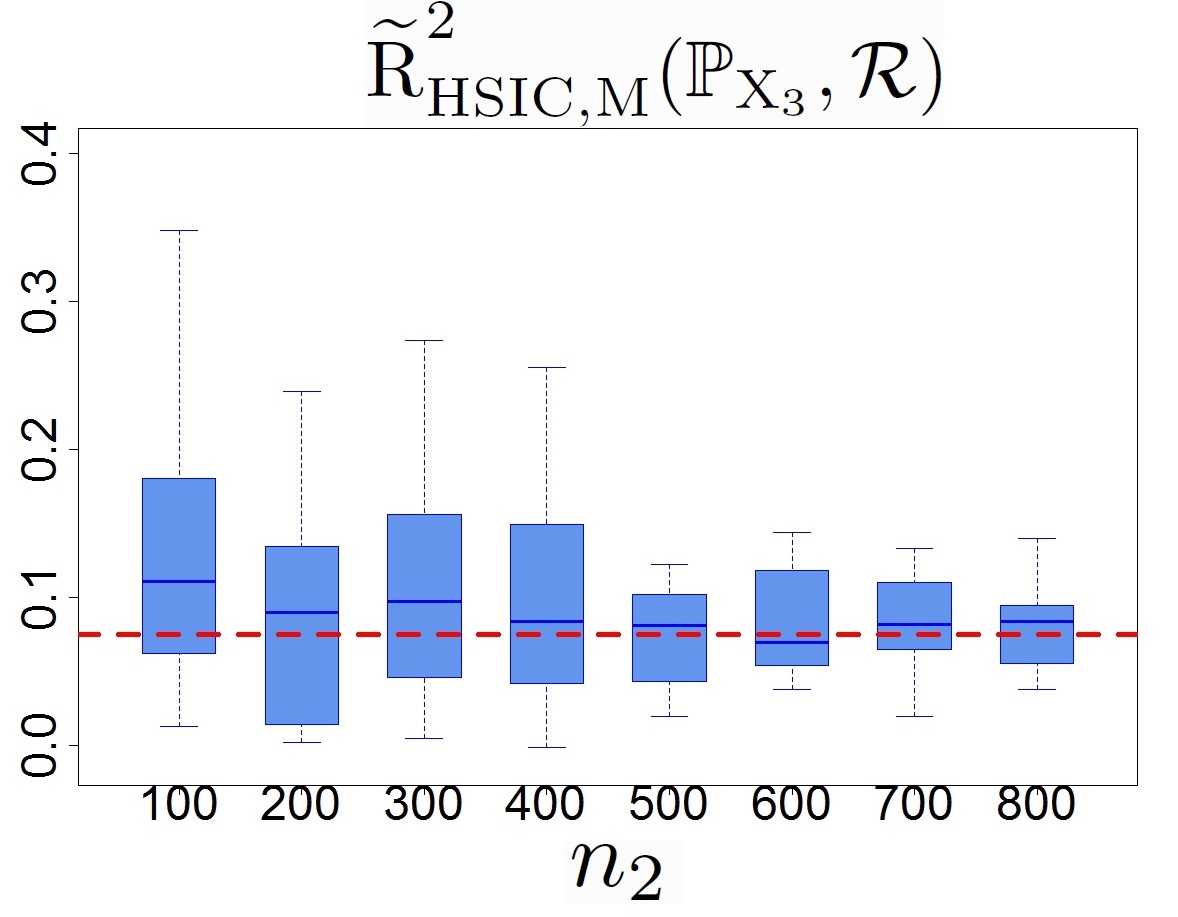}
\caption[]{MACARENa application -- Stability plots of the estimators $\widetilde{\R2}^2_{\HSIC,\M} (\mathbb{P}_{X_k} , \mathcal{R})$, with respect to the sample size $n_2$. Reference values are represented in red dashed lines.}
\label{TestTest}
\end{figure}

\smallskip
Figure \ref{TestTest} presents as a boxplot the mismatch between the values estimated from the initial sample and the  ones estimated from subsamples. We first observe a robustness of estimation: the means of estimators seem to match the value given by the initial sample. We notice also high dispersions for small and medium sizes, (i.e. $n_2 \leq 400$) and small dispersions for medium and big sizes (i.e. $n_2 \geq 500$). Therefore, we conclude that the estimations of GSA2 indices with the sample of $n_2 = 1000$ simulations are consistent, the stability of estimations being satisfactory from $n_2 =700$.

\smallskip
We also check the estimation consistency in terms of input distributions ranking. Table \ref{tauxbonMACA} gives for each subsample size, the rate of times that the ranking matches the one obtained with the initial sample. The results confirm the conclusions drawn from the stability plots. 

\begin{table}[!h]
\begin{center}
\begin{tabular}{|D{2.5cm}|D{1.5cm}|D{1.5cm}|D{1.5cm}|D{1.5cm}|D{1.5cm}|D{1.5cm}|D{1.5cm}|}\hline
$n_2 = 100$ & $n_2 = 200$ & $n_2 = 300$ & $n_2 = 400$ & $n_2 = 500$ &$n_2 = 600$ & $n_2 \geq 700$ \tabularnewline\hline \hline
$45 \%$ & $55\%$ & $70\%$& $75\%$ & $95\%$ & $95\%$ & $100 \%$\tabularnewline \hline
\end{tabular} 
\caption{MACARENa application -- Good ranking rates of estimators $\widetilde{\R2}^2_{\HSIC,\M} (\mathbb{P}_{X_k} , \mathcal{R})$, with respect to the size $n_2$ of the unique sample.}
\label{tauxbonMACA}
\end{center}
\end{table}

\section{Conclusion and Prospect}

In this article, we proposed a new methodology for second-level Global Sensitivity Analysis (GSA2) based on Hilbert-Schmidt Independence Criterion (HSIC). For this, we first proposed new weighted estimators for HSIC, using an alternative sample generated according to a probability distribution which is not the prior distribution of the inputs. We also demonstrated the properties of these new estimators (bias, variance and asymptotic law), which are similar to those of classical estimators. Moreover, their convergence has been illustrated on an analytical example which has also highlighted their ability to correctly rank variables (even for small and medium sample sizes). Subsequently, 2$^{\text{nd}}$-level GSA based on HSIC measures is discussed. When input distributions are uncertain, GSA2 purpose is to assess the impact of these uncertainties on GSA results. In order to perform GSA2, we presented a new \enquote{single-loop} Monte Carlo methodology to address problems raised by GSA2:  characterization of GSA results, definition of 2$^{\text{nd}}$-level HSIC measures and limitation of the calculation budget. This methodology is based on a single sample generated according to a \enquote{reference distribution} (related to the set of all possible distributions). Three options have been proposed for this distribution: mixture law and barycentric laws with respect to the Symmetrical Kullback-Leibler distance or Wasserstein distance. The estimation of  2$^{\text{nd}}$-level HSIC seems to be more accurate using the two first options rather than the Wasserstein barycenter. We also illustrated the great interest of the \enquote{single-loop} approach compared to the \enquote{double-loop} approach. Finally, the whole methodology has been applied to a nuclear test case simulating a severe reactor accident and has shown how GSA2 can provide additional information to classical GSA.   

\smallskip
Several points of the methodology could be more investigated in future research. First, we could focus on comparing Space Filling Design \cite{pronzato2012design, cioppa2002efficient, wang2007review} techniques and Monte Carlo methods for the sampling of input distribution in the case of probabilistic densities (pdf) with uncertain parameters. Indeed, sampling  the uncertain parameters of pdf following a space-filling design could improve the accuracy of the estimators of GSA2 indices. Another interesting perspective would be to build independence tests based on 2$^{\text{nd}}$-level HSIC measures estimators. This could be achieved by identifying the asymptotic distributions of these estimators under the assumption of independence between distributions and GSA1 results.

\smallskip
Furthermore, this new approach for GSA2 could also be compared to the classical approach of epistemic GSA in the framework of Dempster-Shafer theory  \cite{smets1994dempster,alvarez2009reduction}. Indeed, Dempster-Shafer theory gives a description of random variables with epistemic uncertainty, which is to associate with an epistemic variable $Z$ on a set $A$, a mass function representing a probability measure on the set $\mathcal{P} (A)$ of all $A$-subsets. This lack of knowledge is reflected in Dempster-Shafer theory by an upper and lower bound of the cumulative distribution function and can be viewed as 2$^{\text{nd}}$-level of uncertainty. 

\smallskip
An other potential prospect could be to make the connection between our approach and Perturbed-Law based Indices (PLI) \cite{lemaitre2015density,sueur2017sensitivity}. These indices are used to quantify the impact of a perturbation of an input density on the failure probability (probability that a model output exceeds a given threshold). To compare our GSA2 indices with PLI, the probability of failure could be considered as the quantity of interest characterizing GSA results in our methodology. Last but not least, GSA2 method can be compared to the approach proposed in \cite{chabridon2018reliability} which models 2$^{\text{nd}}$-level uncertainties as a uni-level uncertainty on the vector $(\Theta, X)$, where $\Theta$ is the vector of uncertain parameters.   

\section*{Acknowledgments}
We are grateful to Sébastien Da Veiga for his useful ideas and constructive conversations. We also thank Jean-Baptiste Droin for his assistance on the use of MACARENa and Hugo Raguet for his helpful discussions all along this~work.   

\bibliographystyle{plain}
\bibliography{references}

\appendix

\section{Proof of Proposition 1}
\label{sec:AnnexA}

We prove here that
\begin{equation*}
\widetilde{\HSIC} (X_k,Y) = \displaystyle \frac{1}{n^2} \Tr \left(\W \tL_k \W H_1 \tL H_2 \right).
\end{equation*}
Firstly, we evaluate the matrix $ \W \tL_k \W H_1 \tL H_2$ coefficients before computing its trace. The matrix $\W$ being diagonal, we write for $i,j \in \lbrace 1 , \ldots , n \rbrace$:  
\begin{equation*}
( \W \tL_k \W )_{i,j} = ( \tL_k )_{i,j} \W_{i,i} \W_{j,j}. 
\end{equation*} 
The coefficient of the matrix $\W \tL_k \W H_1$ indexed by $i$ and $j$ can therefore be computed:
\begin{align*}
( \W  \tL_k \W H_1 )_{i,j} &= \displaystyle \sum_{r = 1}^n ( \tL_k )_{i,r} \W_{i,i} \W_{r,r} (H_1)_{r,j} \\ &= \displaystyle \sum_{r = 1}^n ( \tL_k )_{i,r}   \W_{i,i}   \W_{r,r}   (\delta_{r,j} - \displaystyle \frac{1}{n} \W_{j,j}) \\ &=  ( \tL_k )_{i,j}   \W_{i,i}    \W_{j,j}  - \displaystyle \frac{1}{n} \sum_{r = 1}^n  ( \tL_k )_{i,r}   \W_{i,i}   \W_{r,r}  \W_{j,j}. 
\end{align*}
Subsequently, the matrix $\W \tL_k \W H_1 \tL$ coefficients are obtained:  
\begin{align*}
( \W \tL_k \W H_1 \tL)_{i,j} &= \displaystyle \sum_{r = 1}^n ( \W \tL_k \W H_1 )_{i,r} \tL_{r,j} \\ &= \displaystyle \sum_{r = 1}^n \left( ( \tL_k )_{i,r}   \W_{i,i}    \W_{r,r} - \displaystyle \frac{1}{n} \sum_{s = 1}^n ( \tL_k )_{i,s}   \W_{i,i}    \W_{s,s}   \W_{r,r} \right) \tL_{r,j} \\ &= \displaystyle \sum_{r = 1}^n  ( \tL_k )_{i,r}   \tL_{r,j}   \W_{i,i}   \W_{r,r}  - \displaystyle \frac{1}{n} \sum_{s = 1}^n  ( \tL_k )_{i,s}   \W_{i,i}   \W_{s,s} \sum_{r = 1}^n \tL_{r,j}   \W_{r,r}. 
\end{align*}
Finally, 
\begin{align*}
( \W \tL_k \W H_1 \tL H_2)_{i,j} &= \displaystyle \sum_{r = 1}^n ( \W \tL_k \W H_1 \tL )_{i,r} (H_2)_{r,j} \\ &= \displaystyle \sum_{r = 1}^n ( \W \tL_k \W H_1 \tL )_{i,r} (\delta_{r,j} - \displaystyle \frac{1}{n} \W_{r,r} ) \\ &= ( \W \tL_k \W H_1 \tL )_{i,j} - \displaystyle \frac{1}{n} \displaystyle \sum_{r = 1}^n ( \W \tL_k \W H_1 \tL)_{i,r} \W_{r,r} \\ &= \displaystyle \sum_{r = 1}^n  ( \tL_k )_{i,r}   \tL_{r,j}   \W_{i,i}   \W_{r,r} - \displaystyle \frac{1}{n} \sum_{1 \leq r,s \leq n} ( \tL_k )_{i,s}  \tL_{r,j}   \W_{i,i}   \W_{s,s}   \W_{r,r} \\ &- \displaystyle \frac{1}{n} \sum_{r = 1}^n \left( \displaystyle \sum_{s = 1}^n  ( \tL_k )_{i,s}   \tL_{s,r}   \W_{i,i}   \W_{s,s} - \displaystyle \frac{1}{n} \sum_{1 \leq p,q \leq n} ( \tL_k )_{i,q}  \tL_{p,r} \W_{i,i}   \W_{q,q}   \W_{p,p} \right) \W_{r,r} \\ &= \displaystyle \sum_{r = 1}^n  ( \tL_k )_{i,r}   \tL_{r,j}   \W_{i,i}   \W_{r,r} - \displaystyle \frac{1}{n} \sum_{1 \leq r,s \leq n} ( \tL_k )_{i,s}  \tL_{r,j}   \W_{i,i}   \W_{s,s}   \W_{r,r} \\ &- \displaystyle \frac{1}{n} \sum_{1 \leq r,s \leq n} ( \tL_k )_{i,s}   \tL_{s,r}   \W_{i,i}   \W_{s,s}   \W_{r,r} + \displaystyle \frac{1}{n^2} \sum_{1 \leq r,p,q \leq n} ( \tL_k )_{i,q}  \tL_{p,r}   \W_{i,i}   \W_{q,q}   \W_{p,p}   \W_{r,r}.
\end{align*}
Summing up the matrix $\W \tL_k \W H_1 \tL H_2$ diagonal terms, then dividing by $n^2$ gives:
\begin{align*}
\nonumber \displaystyle \frac{1}{n^2} \Tr \left( \W \tL_k \W H_1 \tL H_2 \right) &= \displaystyle \frac{1}{n^2} \sum_{1 \leq i,r \leq n}  ( \tL_k )_{i,r}   \tL_{i,r}   \W_{i,i}   \W_{r,r} + \displaystyle \frac{1}{n^4} \sum_{1 \leq i,q \leq n} ( \tL_k )_{i,q}    \W_{i,i}   \W_{q,q} \sum_{1 \leq p,r \leq n} \tL_{p,r}   \W_{p,p}   \W_{r,r} \\ &- \displaystyle \frac{2}{n^3} \sum_{1 \leq i,r,s \leq n} ( \tL_k )_{i,s}   \tL_{i,r}   \W_{i,i}   \W_{s,s}   \W_{r,r}. 
\label{tranx}
\end{align*}
By definition of $\tL_k$, $\tL$ and $\W$, the three terms of the last equation are respectively the estimators defined in Formula \eqref{esHSICmdf}.

\end{document}